\documentclass[12pt]{amsart}
\usepackage{amsmath,latexsym,amssymb,verbatim,euscript}
\theoremstyle{definition}
\newtheorem{thm}{Theorem}
\newtheorem*{thm*}{Theorem}
\newtheorem{prop}{Proposition}[section]

\newtheorem{cor}[prop]{Corollary}
\newtheorem{defn}[prop]{Definition}

\newcommand{\I}{{\mathcal I}}

\newcommand{\Lie}{{\mathcal L}}
\newcommand{\LL}{{\mathbf L}}

\newcommand{\R}{{\mathbb R}}
\newcommand{\ww}{{\mathsf w}}
\newcommand{\vv}{{\mathsf v}}
\newcommand{\vvhat}{{\hat{\vv}}}
\newcommand{\xx}{{\mathbf x}}
\newcommand{\fso}{{\mathfrak so}}
\newcommand{\fs}{{\mathfrak S}}
\newcommand{\fc}{{\mathfrak C}}

\newcommand{\lam}{\lambda}
\newcommand{\w}{\omega}
\newcommand{\ebar}{\bar{e}}
\newcommand{\etabar}{\overline{\eta}}
\newcommand{\thetabar}{\overline{\theta}}
\newcommand{\wbar}{\bar{\omega}}
\newcommand{\zetabar}{\bar{\zeta}}
\newcommand{\Abar}{\bar{A}}
\newcommand{\Dbar}{\bar{D}}
\newcommand{\Fbar}{\overline{\F}}

\newcommand{\pbar}{\bar{p}}
\newcommand{\qbar}{\bar{q}}
\newcommand{\Sbar}{\bar{S}}

\newcommand{\vvbar}{\bar{\vv}}
\newcommand{\xxbar}{\bar{\xx}}
\newcommand{\xbar}{\bar{x}}
\newcommand{\ybar}{\bar{y}}

\newcommand{\fbar}{\bar{f}}
\newcommand{\zbar}{\bar{z}}
\newcommand{\wt}{\widetilde{\w}}
\def\&{\wedge}
\newcommand{\di}{\partial}
\newcommand{\dib}[1]{\dfrac{\di}{\di{#1}}}
\newcommand{\transpose}[1]{{}^t\hskip-2pt{#1}}
\newcommand{\restr}{\negthickspace \mid}
\newcommand{\by}{\times}

\newcommand{\mood}{\mod}

\def\BB/{B\"acklund}
\def\BT/{\BB/ transformation}
\def\MA/{Monge-Amp\`ere}
\def\intprod{\mathbin{\raisebox{.4ex}{\hbox{\vrule height .5pt width
5pt depth 0pt %
           \vrule height 3pt width .5pt depth 0pt}}}}

\newcommand{\bel}[2]{\begin{equation}\label{#1}{#2}\end{equation}}

\newcommand{\om}{\omega}

\newcommand{\calC}{{\mathcal C}}

\newcommand{\calI}{{\mathcal I}}
\newcommand{\calJ}{{\mathcal J}}
\newcommand{\calK}{{\mathcal K}}

\newcommand{\calR}{{\mathcal R}}
\newcommand{\calS}{{\mathcal S}}
\newcommand{\scB}{{\EuScript B}}

\newcommand{\scG}{{\EuScript G}}
\newcommand{\scF}{{\EuScript F}}
\newcommand{\scM}{{\EuScript M}}
\newcommand{\scN}{{\EuScript N}}
\newcommand{\scP}{{\EuScript P}}
\newcommand{\scrU}{{\EuScript U}}
\newcommand{\F}{{\scF}}
\newcommand{\bb}{\mathbb}
\newcommand{\B}{\scB}
\newcommand{\Blam}{\B_{\lam_0}}
\newcommand{\G}{\scG}

\newenvironment{ex}{\begin{trivlist} \item[] {\it Example:} }{
\end{trivlist}  }
\newenvironment{contex}{\begin{trivlist} \item[] {\it Example
(cont'd):} }{ \end{trivlist}  }
\newenvironment{rem}{\begin{trivlist} \item[] {\it Remark.}}{\end{trivlist}}
\newenvironment{pf}{\begin{trivlist} \item[] {\it Proof.} }{\qed \end{trivlist} } 

\begin{document}
\title{Parametric B\"acklund Transformations I: Phenomenology}
\author{Jeanne N. Clelland}
\address{Dept. of Mathematics, 395 UCB, University of
Colorado\\
Boulder, CO 80309-0395}
\email{Jeanne.Clelland@colorado.edu}
\author{Thomas A. Ivey}
\address{Dept. of Mathematics, College of Charleston\\
66 George St., Charleston SC 29424-0001}
\email{IveyT@cofc.edu}


\begin{abstract}
We begin an exploration of parametric B\"acklund transformations for
hyperbolic Monge-Amp\`ere systems.  We compute invariants for
such transformations and explore the behavior of four examples
regarding their invariants, symmetries, and conservation laws.  We
prove some preliminary results and indicate directions for further
research.
\end{abstract}

\maketitle

\section{Introduction}\label{intro}
In this paper we will study parametric \BB/ transformations between
second-order hyperbolic {\it Monge-Amp\`ere equations} for one
function of two variables.
These are partial differential equations (PDE) of the form
\[ A (z_{xx}z_{yy} - z_{xy}^2) + B z_{xx} + 2C z_{xy} + D z_{yy} +
E = 0 \]
where the coefficients $A,B,C,D,E$ are functions of the variables
$x,y,z,z_x, z_y$.  The equation is {\em hyperbolic} if it has
distinct, real characteristics at each point, i.e., if $AE - BD +
C^2 > 0$.

Before giving a formal definition, we consider the following
classical example.

\begin{ex} Consider the pair of PDEs
\begin{align}
z_x - \zbar_x & = \lam\, \sin(z + \zbar) \label{SGEback} \\
z_y + \zbar_y & = \frac{1}{\lam}\, \sin(z - \zbar) \notag
\end{align}
for functions $z(x,y),\ \zbar(x,y)$, with $\lambda$ a
nonzero real number.  By differentiating the first equation with respect
to $y$ and the second equation with respect to $x$ and then adding or
subtracting the two equations, it can be shown that
any smooth functions $z(x,y),\ \zbar(x,y)$ satisfying these equations
must both be solutions of the sine-Gordon equation
\begin{equation}
    z_{xy} = \tfrac{1}{2}\sin (2z) . \label{SGE}
\end{equation}
   From another point of view, if
$z(x,y)$ is a known solution of \eqref{SGE},  then \eqref{SGEback} is a
compatible, overdetermined system for $\zbar(x,y)$.
Solving this system, which only requires integrating
{\em ordinary} differential equations, yields new solutions of \eqref{SGEback}
which are {\it \BB/ transformations} of the original solution.
For example, taking $z(x,y) = 0$ yields
\[ \zbar(x,y) = \tan^{-1}(e^{-(\lam x + \tfrac{1}{\lam}y + c)}). \]
These are the {\em 1-soliton} solutions of the sine-Gordon equation.
\end{ex}

The sine-Gordon equation is a classical example of an integrable
system.  The remarkable features typical of such systems include
an infinite sequence of conservation laws, a method of solution
by inverse scattering \cite{DJ}, and \BB/ transformations.
(Indeed, the most interesting examples of \BB/ transformations come from
integrable systems.)
Note that the arbitrary parameter $\lam$ appearing in the \BB/
transformation \eqref{SGEback} does not appear in the sine-Gordon
equation itself.  This parameter plays a critical role in the theory of
\BB/ transformations of integrable systems; for instance, it is
possible to derive the infinite hierarchy of conservation laws for
the sine-Gordon equation by performing a series expansion of the equations
\eqref{SGEback} in terms of $\lam$ \cite{WSK}.  Furthermore, \eqref{SGEback}
is equivalent to the sine-Gordon AKNS system (with $\lam$ as spectral
parameter), which is
the centerpiece of solution by inverse scattering \cite{rogers}.

\begin{rem} The system \eqref{SGEback} is technically an {\em
auto}-\BB/ transformation for
sine-Gordon.  In general, \BB/ transformations can take solutions
of one PDE and produce solutions to a different PDE.  For example,
suppose that   $z(x,t)$
is a solution of the nonlinear heat equation
$z_t=(f(z) z_x)_x.$
If $y(x,t)$ is a solution to the equations
$$y_x =z,\qquad y_t = f(z) z_x$$
and we set $\zbar=1/z$, then $\zbar$ is a solution to another
nonlinear heat equation
$\zbar_t=(\fbar(\zbar) \zbar_y)_y$, where $\fbar$ is related to $f$
by $z f(z)=\zbar \fbar(\zbar)$.  This example also shows that, in general,
\BB/ transformations may involve a change of independent variable.
\end{rem}

As in \cite{C01b}, we will generalize the notion of a Monge-Amp\`ere
equation slightly, considering \BB/ transformations
between certain exterior differential systems known as Monge-Amp\`ere systems.

An {\em exterior differential system} (EDS) on a manifold $\scM$ is a
differentially
closed ideal ${\calI}$ in the graded algebra, under wedge product, of
differential forms on $\scM$. Any system of partial differential
equations can be
formulated as an EDS $\calI$, and solutions of the PDE system
correspond to {\em
integral manifolds} of $\calI$, i.e., submanifolds $\scN \subset
\scM$ which satisfy
the condition that all the forms in $\calI$ vanish when pulled back to $\scN$.

\begin{defn}\label{MAdefn}
A {\em Monge-Amp\`ere system} $\calI$ is an EDS
on a 5-dimensional manifold $\scM$ such that $\calI$ is locally generated by a
contact form $\theta$ (i.e., a 1-form $\theta$ with the property that
$\theta \& d\theta \& d\theta \neq 0$), the 2-form $\Theta = d\theta$,
and another 2-form $\Psi$ which is linearly independent from $\Theta$
and wedge products involving $\theta$.
\end{defn}

A Monge-Amp\`ere system $\calI$ is {\em hyperbolic} if
the homogeneous quadratic equation
$$(\mu\, \Theta + \nu\, \Psi) \& (\mu\, \Theta + \nu\, \Psi)
\equiv 0 \mod{\theta}$$
has linearly independent real roots $[\mu_i,\nu_i]$, $i=1,2$.
This condition agrees with the traditional
definition of hyperbolicity, and it implies that there are two
independent linear combinations
$\mu_i\, \Theta + \nu_i\, \Psi +
\beta_i \wedge \theta$
which are {\em decomposable}---i.e., they can each be written as
a wedge product of a pair of 1-forms.  Associated to $\calI$ are two
{\em Monge characteristic systems}, each containing the factors of one
these decomposable 2-forms, together with $\theta$.   (The
characteristic systems must each contain $\theta$ because the factors
of the decomposable 2-forms are only well-defined modulo $\theta$.)

\begin{rem}  Those EDS which are generated algebraically by 1-forms and their
exterior derivatives are known as {\it Pfaffian systems}; the dimension of
the space of 1-forms in the system is the {\em rank} of the system.
Any hyperbolic \MA/ system is equivalent, under prolongation, to a type of
rank three Pfaffian system known as a hyperbolic exterior
differential system \cite{BGH95}.
\end{rem}

\begin{contex}
The sine-Gordon equation \eqref{SGE} may be described
as a hyperbolic Monge-Amp\`ere system $\calI$ on $\mathbb{R}^5$ (with
coordinates $(x,y,z,p,q)$)
generated by the forms
\begin{gather}
    \theta = dz - p\,dx - q\, dy \notag \\
    \Theta = -dp \& dx - dq \& dy \\
    \Psi  = [dp - \tfrac{1}{2} \sin(2z)\, dy] \& dx \notag.
\label{SGEMA}
\end{gather}
Two-dimensional integral manifolds of this system that satisfy the
independence condition
$dx \& dy \neq 0$ are naturally in one-to-one correspondence with solutions of
\eqref{SGE}.

Note that $\Psi$ is decomposable; another decomposable linear
combination of $\Psi$ and $\Theta$ is $\Psi + \Theta = -[dq - \tfrac{1}{2}
\sin(2z)\, dx] \& dy$.  Thus, the two characteristic systems are
\begin{align*}
\calC_1&=\{dz-p\,dx-q\,dy,\ dp - \tfrac{1}{2} \sin(2z)\, dy,\ dx\},\\
\calC_2&=\{dz-p\,dx-q\,dy,\ dq - \tfrac{1}{2} \sin(2z)\, dx,\ dy\}.
\end{align*}
\end{contex}

\bigskip
Given two hyperbolic Monge-Amp\`ere systems $(\scM_1, \calI_1)$ and $(\scM_2,
\calI_2)$, we define a (non-parametric) {\em \BB/ transformation}
between $(\scM_1,\calI_1)$ and $(\scM_2, \calI_2)$ to be a
6-dimensional submanifold $\scB \subset \scM_1 \times \scM_2$
which has the following properties:
\begin{list}{\arabic{enumi}.}
{\usecounter{enumi}}
\item{The natural projections $\pi_1:\scB \rightarrow \scM_1$ and
$\pi_2: \scB \rightarrow \scM_2$ are submersions.

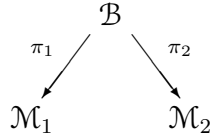
\setlength{\unitlength}{2pt}
\begin{center}
\begin{picture}(40,30)(0,0)
\put(5,5){\makebox(0,0){$\scM_1$}}
\put(35,5){\makebox(0,0){$\scM_2$}}
\put(20,25){\makebox(0,0){$\scB$}}
\put(16,21){\vector(-3,-4){9}}
\put(24,21){\vector(3,-4){9}}
\put(7,18){\makebox(0,0){$\scriptstyle{\pi_1}$}}
\put(33,18){\makebox(0,0){$\scriptstyle{\pi_2}$}}
\end{picture}
\end{center}
}
\item{The pullbacks to $\scB$ of the forms $\Theta_1, \Theta_2,
\Psi_1, \Psi_2$ satisfy the condition that
\begin{equation}\label{Jcond}
    \{ \Psi_1, \Psi_2 \} \equiv \{ \Theta_1, \Theta_2 \}
      \mod{\theta_1, \theta_2}.
\end{equation}
Since $\Theta_1, \Psi_1$ are linearly independent forms (as are
$\Theta_2, \Psi_2$), this condition implies that
\begin{equation*}
    \{ \Theta_1, \Psi_1 \}  \equiv \{ \Theta_2, \Psi_2  \}
      \mod{\theta_1, \theta_2}.
\end{equation*}
This second equation is really the desired
property; the first equation ensures that, in
addition, the forms $\Theta_1, \Theta_2$ are linearly independent.  }
\end{list}

That this definition captures the properties of the sine-Gordon
example may be seen as follows: suppose that $\scN
\hookrightarrow \scM_1$ is a two-dimensional integral manifold of $\calI_1$.
The inverse image
$\pi_1^{-1}(\scN)$ is a three-dimensional submanifold of $\scB$.  Now
consider the
restriction of $\pi_2^* (\calI_2)$ to $\pi_1^{-1}(\scN)$.  By
Property (2) above, the restriction of $\pi_2^* (\calI_2)$ to
$\pi_1^{-1}(\scN)$ is a
Frobenius system (i.e., a Pfaffian system which is
generated {\em algebraically} by its 1-forms).
By the Frobenius Theorem, $\pi_1^{-1}(\scN)$ is foliated by two-dimensional
integral manifolds of $\pi_2^* (\calI_2)$, each of which
projects to an integral manifold of $(\scM_2, \calI_2)$; moreover, these
integral manifolds can be constructed by solving ODEs.

\begin{contex} Let $\scM_1 = \scM_2 = \bb{R}^5$ with coordinates
$(x,y,z,p,q)$ on $\scM_1$ and $(\bar{x},\bar{y},\bar{z},\bar{p},\bar{q})$
on $\scM_2$.  Let $\calI_1$ and $\calI_2$ be copies of the sine-Gordon system
$\calI$ described above, defined on $\scM_1$ and $\scM_2$, respectively.
For any fixed value of $\lam$,
the \BB/ transformation \eqref{SGEback} may be regarded as the
6-dimensional submanifold $\scB \subset \scM_1 \times \scM_2$ defined by
the four equations
\begin{gather*}
\bar{x} = x, \qquad \bar{y} = y, \\
p - \bar{p} = \lam\, \sin(z + \zbar),  \\
q + \bar{q} = \frac{1}{\lam}\, \sin(z - \zbar).
\end{gather*}
\end{contex}

In \cite{C01b}, the first author used Cartan's method of equivalence
to classify all homogeneous
non-parametric \BB/ transformations of hyperbolic Monge-Amp\`ere
systems.  (Here, ``homogeneous" means that the $G$-structure constructed
via the method of equivalence is assumed to be acted on by a transitive
symmetry group, implying that all structure functions are constants.)
The restriction to the homogeneous case was necessary in
order to simplify the computations in the method of equivalence so
that the classification could be completed explicitly.  Unfortunately,
homogeneity is a significant restriction.  In particular, the
\BB/ transformation for the sine-Gordon equation is not
homogeneous, so we expect that there are many interesting
transformations which do not appear in that classification.

In the present paper, we wish to explore \BB/ transformations
that depend explicitly on an arbitrary parameter, as in the
sine-Gordon case.  To this end, we define a {\em parametric
\BB/ transformation} between two hyperbolic Monge-Amp\`ere
systems $(\scM_1, \calI_1)$ and $(\scM_2, \calI_2)$ to be a
7-dimensional submanifold
\[ \scP \subset \scM_1 \times \scM_2 \times \bb{R}, \]
with coordinate $\lambda$ on the $\bb{R}$ factor,
such that for each $\lam_0 \in \bb{R}$, the set
\[ \Blam = \scP \cap (\scM_1 \times \scM_2 \times \{\lam_0\}) \]
is a non-parametric \BB/ transformation.  (We will let $\pi_1,\pi_2$
also denote the submersions from $\scP$ to $\scM_1,\scM_2$ respectively.)
   We will
associate a canonical $G$-structure to such a transformation, and we
will see that the transformation's nontrivial dependence on $\lam$ will
be equivalent to the nonvanishing of certain invariants appearing in
the structure equations.

Ultimately, we hope to explore the relationship between \BB/
transformations and other features typically associated with
integrable systems, such as symmetries of the system or an infinite hierarchy
of conservation laws.  In this paper we will explore in some detail
the behavior of four examples of parametric \BB/ transformations:
\begin{itemize}
\item the sine-Gordon transformation described above;
\item a transformation, given by Zvyagin \cite{zvyagin}, between the
wave equation and a quasilinear PDE which we refer to as
Goursat's equation;
\item the classical auto-\BB/ transformation between surfaces of constant Gauss
curvature $K=-1$ in $\bb{E}^3$;
\item an auto-\BB/ transformation between timelike surfaces of constant mean
curvature (CMC) equal to $1$ in the negatively curved Lorentzian space form
$\mathbf{H}^{2,1}$ (also known as {\em anti-de Sitter space}).
\end{itemize}

We will see that:
\begin{itemize}
\item In the first two examples, nontrivial dependence on a parameter arises
from symmetries of the underlying \MA/ systems
$(\scM_1, \calI_1),\, (\scM_2, \calI_2)$, but no such symmetries
lead to the parametric dependence in the other two examples.
\item In the first two examples, each characteristic system contains
a rank one integrable subsystem, whereas in the $K=-1$ example
none of the characteristic systems contains such a subsystem, and
in the CMC example exactly one of the characteristic systems for each of
$\I_1$ and $\I_2$ contains such a subsystem.  For this reason we
regard the CMC example
as somehow exhibiting behavior ``in between'' that of the other examples.
\item In the sine-Gordon example, each \MA/ system has a 3-dimensional space
of conservation laws, while the space of conservation laws for the system
$\calJ = \{\theta_1,\theta_2,\Theta_1,\Theta_2\}$ on $\Blam$ is
4-dimensional, and is the
union of the pullbacks of the spaces of conservation laws from
$\scM_1$ and $\scM_2$.
By contrast, in the $K=-1$ example each \MA/ system has a
6-dimensional space of conservation laws, and the space of conservation laws
on $\Blam$ is a 7-dimensional space in which the union of the pullbacks forms a
6-dimensional subspace.  (In the other two examples the space of conservation
laws for each \MA/ system is infinite-dimensional.)
\end{itemize}

In the sequel to this paper, we plan to determine the relationship
between these properties and certain invariants appearing in the structure
equations of the $G$-structure associated to parametric \BB/ transformations in
general.  Some initial results in this direction will appear in
section \ref{concon}.

\section{A canonical $G$-structure}\label{Gstrsec}
In this section we will associate a canonical $G$-structure to a
parametric \BB/ transformation between two \MA/ systems $(\scM_1,\calI_1),\,
(\scM_2, \calI_2)$.

Suppose that $\scP \subset \scM_1 \times \scM_2 \times \R$ is a parametric \BB/
transformation as described above.
Let $\thetabar_1$, $\thetabar_2$ denote locally defined 1-forms on $\scP$
which are the pullbacks of the contact forms on $\scM_1,\scM_2$, respectively,
let $\zetabar$ be an integrable 1-form on $\scP$ whose pullback
to each slice $\Blam$ vanishes identically, and let
$\calK$ denote the span of $\{\thetabar_1, \thetabar_2, \zetabar\}$.
(For the present, we will use
barred letters to denote forms and functions on $\scP$ and unbarred letters to
denote the corresponding semi-basic forms and functions
defined on a $G$-structure over $\scP$.)

Recall that the {\em Cartan system} (also called the {\em retracting
space} \cite{EDS})
of a Pfaffian system $\I$ is the smallest Frobenius system containing $\I$.  In
particular, the Cartan system of a single 1-form $\theta$ is
spanned by $\theta$ and the forms which appear in $d\theta$
mod $\theta$, and the dimension of this Cartan system is the
{\em Pfaff rank} of the 1-form. Since
\bel{zetacond}{d\zetabar \equiv 0 \mod \zetabar,}
the Pfaff rank of $\zetabar$ is one.
Since each $\thetabar_i$ is the pullback of a contact form on a
5-dimensional manifold, its Cartan system has rank 5 and
consists of 1-forms which are {\em semi-basic} for $\pi_i$,
i.e. they vanish on vectors tangent to the fibers of $\pi_i$.

Let $\calR_1,\calR_2$ denote the Cartan systems
of $\thetabar_1,\thetabar_2$, respectively.
We are free to modify each of $\thetabar_1,\thetabar_2$, and $\zetabar$
by multiplying by a nonzero function, and if we do so
their Pfaff ranks and Cartan systems do not change.

Let $\calC_{i1},\calC_{i2}$ denote the Monge characteristic systems of $\I_i$.
Since each slice $\Blam$ defines a \BT/ between $\I_1$ and $\I_2$, there
are 1-forms $\wbar^1,\wbar^2,\wbar^3,\wbar^4$ on $\scP$ such that
\begin{alignat*}{2}
\pi_1^*\calC_{11} &= \{\thetabar_1,\wbar^1,\wbar^2\},\qquad &
\pi_2^*\calC_{21} &= \{\thetabar_2,\wbar^1,\wbar^2\},\\
\pi_1^*\calC_{12} &= \{\thetabar_1,\wbar^3,\wbar^4\},\qquad &
\pi_2^*\calC_{22} &= \{\thetabar_2,\wbar^3,\wbar^4\},
\end{alignat*}
{\em when restricted to $\B_{\lambda_0}$.}
(See \cite{C01b} for a detailed construction of this coframe on any slice.)
Since $d\theta_i \mod \theta_i$
must be a linear combination of decomposables in $\calC_{i1}$ and
$\calC_{i2}$, we have
\begin{equation}\label{dthetazeta}
\begin{aligned}
d\thetabar_1 &\equiv \Abar_1\, \wbar^1 \& \wbar^2 +
\Abar_3\, \wbar^3 \& \wbar^4\mod \thetabar_1,\zetabar \\
d\thetabar_2 &\equiv \Abar_4\, \wbar^1 \& \wbar^2 + \Abar_2\, \wbar^3 \&
\wbar^4 \mod \thetabar_2,\zetabar\end{aligned}
\end{equation}
for some nonzero functions $\Abar_i$ on $\scP$ such that
$\Abar_1\Abar_2-\Abar_3\Abar_4 \ne 0$ at each point.  This implies the
analogue of \eqref{Jcond},
$$
\{\pi_1^* \bar{\Psi}_1, \pi_2^* \bar{\Psi}_2 \}
\equiv \{d\thetabar_1, d\thetabar_2\} \mood \calK.$$

We are free to modify the $\wbar^i$ by scaling and adding multiples
of $\zetabar$,
and we will use this to arrange that $\Abar_3=1=\Abar_4$ and that, on $\scP$,
\bel{Mcond}{\pi_1^*\calC_{11} = \{\thetabar_1,\wbar^1,\wbar^2\}\qquad
\pi_2^*\calC_{22} = \{\thetabar_2,\wbar^3,\wbar^4\}.}
(This condition determines the spans
$\{\wbar^1,\wbar^2\}$ and $\{\wbar^3,\wbar^4\}$ uniquely.)
Consequently, there are functions $\Dbar_i$ on $\scP$ such that
\bel{thetaAD}{\begin{aligned}
d\thetabar_1 &\equiv \Abar_1\, \wbar^1 \& \wbar^2 + \wbar^3
\& \wbar^4 + \zetabar \& (\Dbar_3\, \wbar^3 + \Dbar_4\, \wbar^4)
\mood \thetabar_1 
\\ d\thetabar_2 &\equiv
  \wbar^1 \& \wbar^2 + \Abar_2\, \wbar^3 \& \wbar^4 +
  \zetabar \& (\Dbar_1\, \wbar^1 + \Dbar_2\, \wbar^2) \mood
\thetabar_2
\end{aligned}}
with $\Abar_1 \Abar_2 \neq 1$, and
\bel{Mothercond}{
\pi_1^*\calC_{12} =\{\thetabar_1,\ \wbar^3+\Dbar_4\zetabar,\
\wbar^4-\Dbar_3\zetabar\},
\qquad
\pi_2^*\calC_{21} = \{\thetabar_2,\ \wbar^1+\Dbar_2\zetabar,\
\wbar^2-\Dbar_1\zetabar\}.
}

Now let $\scG$ be the sub-bundle of the coframe bundle of $\scP$
consisting of all coframes
$(\zetabar,\, \thetabar_1,\, \thetabar_2,\, \wbar^1,$ $\wbar^2,\,
\wbar^3,\, \wbar^4)$
satisfying the conditions \eqref{zetacond}, \eqref{Mcond}, \eqref{thetaAD},
\eqref{Mothercond} above.  Then $\scG$ is
a $G$-structure on $\scP$, where $G\subset GL(7,\R)$ is the
9-dimensional group of
matrices of the form $$\left\{\begin{bmatrix}
\begin{smallmatrix}t&0&0\\0&s&0\\0&0&r\end{smallmatrix}&0&0\\
0& R& 0\\
0& 0& S
\end{bmatrix}: R, S \in GL(2,\R), \ \det(R)=r,\ \det(S)=s,\ t \neq
0\right\}.$$
There are canonical 1-forms
$(\zeta,\, \theta_1,\, \theta_2,\, \w^1,$ $\w^2,\,
\w^3,\, \w^4)$ on $\scG$ defined by the property that for any
section
$\sigma = (\zetabar,\, \thetabar_1,\, \thetabar_2,\, \wbar^1,\,\wbar^2,\,
\wbar^3,\, \wbar^4): \scP \to \scG$
we have
\bel{reprop}{\sigma^*(\zeta) = \zetabar,\ \sigma^*(\theta_i) =
\thetabar_i,\ \sigma^*(\w^i) = \wbar^i.}
These forms are semi-basic for the projection from
$\scG$ to $\scP$ and will henceforth be referred
to as the {\em semi-basic forms} on $\scG$. The
reproducing property \eqref{reprop} implies that
\bel{dadapt}{\begin{aligned}d\theta_1 &\equiv A_1\, \w^1 \& \w^2 + \w^3
\& \w^4 + \zeta \& (D_3\, \w^3 + D_4\, \w^4) \mood \theta_1\\
    d\theta_2 &\equiv \w^1 \& \w^2 + A_2\, \w^3 \& \w^4 +  \zeta \& (D_1\,
    \w^1 + D_2\, \w^2) \mood \theta_2
\end{aligned}}
where $A_i,D_i$ are functions on $\scG$ such that $\sigma^*(A_i)=\Abar_i$
and $\sigma^*(D_i)=\Dbar_i$ for any section $\sigma$.
Using $\calR_1,\calR_2$ to denote the Cartan systems of $\theta_1,\theta_2$
on $\scG$, we have
$$\calR_1=\{ \theta_1, \w^1, \w^2, \wt^3, \wt^4 \}, \qquad
\calR_2 = \{ \theta_2, \wt^1, \wt^2, \w^3,\w^4\},$$
where, for the sake of convenience, we define
$$\wt^1 = \w^1 + D_2\, \zeta, \quad
  \wt^2 = \w^2 - D_1\, \zeta,\quad \wt^3= \w^3 + D_4\, \zeta,\quad \wt^4 =
\w^4 - D_3\, \zeta.$$

We will now proceed with the method of equivalence by trying
to normalize terms appearing in the exterior derivatives of the
semi-basic forms on $\scG$.
First, we must compute the exterior derivatives of the $\w^i$.

Note that the pullbacks of the Monge characteristic systems of $\I_1$ are
given by $\{\theta_1,\w^1,\w^2\}$ and $\{\theta_1,\wt^3,\wt^4\}$,
while the pullbacks of those of $\I_2$ are
$\{\theta_2,\wt^1,\wt^2\}$ and $\{\theta_2,\w^3,\w^4\}$.  This implies that
\bel{nontcond}{
\begin{aligned}d\w^1\equiv d\w^2&\equiv 0 \mod
\theta_1,\w^1,\w^2,\Lambda^2(\calR_1)\\
d\w^3\equiv d\w^4&\equiv 0 \mod \theta_2,\w^3,\w^4,\Lambda^2(\calR_2)
\end{aligned}}
and
\bel{tcond}{
\begin{aligned}d\wt^1\equiv d\wt^2&\equiv 0 \mod
\theta_2,\wt^1,\wt^2,\Lambda^2(\calR_2)\\
d\wt^3\equiv d\wt^4&\equiv 0 \mod \theta_1,\wt^3,\wt^4,\Lambda^2(\calR_1).
\end{aligned}}
Because $\zeta$ is integrable, \eqref{tcond} implies that
\bel{nontcondzeta}{
\begin{aligned}d\w^1\equiv d\w^2&\equiv 0 \mod
\zeta,\theta_2,\w^1,\w^2,\Lambda^2(\calR_2)\\
d\w^3\equiv d\w^4&\equiv 0 \mod \zeta,\theta_1,\w^3,\w^4,\Lambda^2(\calR_1).
\end{aligned}}
Together, \eqref{nontcond} and \eqref{nontcondzeta} imply that
\begin{align*}
& \left. \begin{array}{l}
d\w^1 \equiv B_1\, \theta_1 \& \theta_2 + C_1\, \wt^3 \& \wt^4  +
E_1\, \zeta \& \theta_1 \\[0.1in]
d\w^2 \equiv B_2\, \theta_1 \& \theta_2 + C_2\, \wt^3 \& \wt^4  +
E_2\, \zeta \& \theta_1
\end{array} \right\} \mod{\om^1, \om^2} \\[0.1in]
& \left. \begin{array}{l}
d\w^3 \equiv B_3\, \theta_1 \& \theta_2 + C_3\, \wt^1 \& \wt^2 +
E_3\, \zeta \& \theta_2 \\[0.1in]
d\w^4 \equiv B_4\, \theta_1 \& \theta_2 + C_4\, \wt^1 \& \wt^2 +
E_4\, \zeta \& \theta_2
\end{array} \right\} \mod{\om^3, \om^4}
\end{align*}
for some functions $B_i,C_i$ and $E_i$ on $\scG$.

\newcommand{\cofcolumn}{{\begin{bmatrix} \zeta \\ \theta_1
\\ \theta_2 \\ \w^1 \\ \w^2 \\ \w^3 \\ \w^4 \end{bmatrix}}}
A standard argument shows that on
$\scG$ there exist 1-forms $\alpha_i,\,\beta_i,\, \gamma$
(referred to as {\em connection forms}), linearly independent from
the semi-basic forms,
such that
\begin{equation}
d\cofcolumn
= - \begin{bmatrix} \gamma & 0 & 0 & 0 & 0 & 0 & 0 \\
                       0      & \beta_1+\beta_4 & 0 & 0 & 0 & 0 & 0 \\
                       0 & 0 & \alpha_1+\alpha_4& 0 & 0 & 0 & 0 \\
                       0 & 0 & 0 & \alpha_1 & \alpha_2 & 0 & 0 \\
                       0 & 0 & 0 & \alpha_3 & \alpha_4 & 0 & 0 \\
                       0 & 0 & 0 & 0 & 0 & \beta_1 & \beta_2 \\
0 & 0 & 0 & 0 & 0 & \beta_3 & \beta_4 \end{bmatrix} \& \cofcolumn +
\begin{bmatrix} \Upsilon \\ \Theta_1 \\ \Theta_2 \\ \Omega^1 \\
\Omega^2 \\ \Omega^3 \\ \Omega^4 \end{bmatrix}. \label{bigstreq}
\end{equation}
These are the {\em structure equations} of $\scG$.  The column on the far right
consists of 2-forms which are sums of wedge products of
pairs of semi-basic forms on $\scG$; these forms are collectively known as
{\em torsion.}  Our assumptions about the $G$-structure limit the kinds of
terms that can occur as torsion, and we may also absorb some torsion
terms by adding semi-basic forms to the connection forms.  In particular,
we can arrange that
\begin{align*}
\Upsilon &= 0 \\
\Theta_1 &= F_1\, \zeta \&\theta_1 + A_1\, (\w^1 -
C_1\, \theta_1) \& (\w^2 - C_2\, \theta_1) + (\w^3 + D_4\, \zeta)
\&(\w^4 - D_3\, \zeta) \\
\Theta_2 &= F_2\, \zeta \& \theta_2 + A_2\, (\w^3 -
C_3\, \theta_2) \& (\w^4 - C_4\, \theta_2) + (\w^1 +
D_2\,\zeta)\&(\w^2 - D_1\, \zeta)\\
\Omega_1 &= B_1\, \theta_1 \& \theta_2 + C_1\, (\w^3 + D_4\, \zeta)
\&(\w^4 - D_3\, \zeta)
   + E_1\, \zeta \& \theta_1  \\
\Omega_2 &= B_2\, \theta_1 \& \theta_2 + C_2\, (\w^3 + D_4\, \zeta)
\&(\w^4 - D_3\, \zeta)
+ E_2\, \zeta \&\theta_1 \\
\Omega_3 &= B_3\, \theta_1 \& \theta_2 + C_3\, (\w^1 + D_2\,
\zeta)\&(\w^2 - D_1\, \zeta)
+ E_3\, \zeta \&\theta_2 \\
\Omega_4 &= B_4\, \theta_1 \& \theta_2 + C_4\, (\w^1 + D_2\,
\zeta)\&(\w^2 - D_1\, \zeta)
+ E_4\, \zeta \& \theta_2
\end{align*}
for some functions $A_i, B_i, C_i, D_i, E_i, F_i$ on $\scG$.  (As
noted previously, we must have $A_1,\, A_2$ nonzero and $A_1 A_2 \neq
1$.)

\newcommand{\bonevec}{\left[\begin{matrix}{B_1\\
B_2}\end{matrix}\right]}
\newcommand{\colvec}[2]{\left[\begin{matrix}{#1 &
#2}\end{matrix}\right]}

The right action of $G$ on sections of $\scG$ is $g\cdot
\sigma=g^{-1} \sigma$.
This induces an action on the torsion coefficients, as follows.
The functions $A_1,\,A_2,\,F_1,\,F_2$ are acted on by scaling:
$$ A_1 \to rs^{-1} A_1, \qquad A_2 \to r^{-1}s A_2,
\qquad F_1 \to t F_1,  \qquad F_2 \to t F_2.$$
(In particular, the product $A_1 A_2$ is invariant on the fibers of $\scG$,
so the condition $A_1 A_2 \neq 1$ makes sense.)
The remaining coefficients occur naturally in pairs
as components of vectors and are acted on as follows:
\begin{alignat*}{4}
\begin{bmatrix} B_1 \\ B_2 \end{bmatrix} & \to &
rs R^{-1} &\begin{bmatrix} B_1 \\ B_2 \end{bmatrix} \qquad \qquad &
\begin{bmatrix} B_3 \\ B_4 \end{bmatrix} & \to &
rs S^{-1} &\begin{bmatrix} B_3 \\ B_4 \end{bmatrix}  \\[0.1in]
\begin{bmatrix} C_1 \\ C_2 \end{bmatrix} & \to &
s R^{-1} &\begin{bmatrix} C_1 \\ C_2 \end{bmatrix}  &
\begin{bmatrix} C_3 \\ C_4 \end{bmatrix} & \to &
r S^{-1} &\begin{bmatrix} C_3 \\ C_4 \end{bmatrix}  \\[0.1in]
\begin{bmatrix} D_1 \\ D_2 \end{bmatrix} & \to &
(r^{-1}t)\,\transpose{R} &\begin{bmatrix} D_1 \\ D_2 \end{bmatrix} &
\begin{bmatrix} D_3 \\ D_4 \end{bmatrix} & \to &
(s^{-1}t)\,\transpose{S} &\begin{bmatrix} D_3 \\ D_4 \end{bmatrix}  \\[0.1in]
\begin{bmatrix} E_1 \\ E_2 \end{bmatrix} & \to &
st R^{-1} &\begin{bmatrix} E_1 \\ E_2 \end{bmatrix} &
\begin{bmatrix} E_3 \\ E_4 \end{bmatrix} & \to &
rt S^{-1} &\begin{bmatrix} E_3 \\ E_4 \end{bmatrix}.
\end{alignat*}
In particular, we observe that the vanishing of any of these vectors
is a well-defined condition on $\scP$, as is the linear dependence or
independence of certain pairs of vectors, e.g.
$[B_1\  B_2]$ and $[C_1 \ C_2]$ or $[E_1 \ E_2]$.
As well, the dot product of $[D_1 \ D_2]$ with any one of these
three is a relative invariant.

We can make some elementary observations regarding how the vanishing
of some of these vectorial invariants implies the vanishing of others.
\begin{prop} If the vector $[C_1\  C_2]$ vanishes identically on
$\scP$,
then so do $[B_1\  B_2]$ and $[E_1\  E_2]$.  Similarly, if $[C_3\  C_4]$
vanishes, then so do $[B_3\  B_4]$ and $[E_3\  E_4]$.
\end{prop}
\begin{pf}
Suppose that $C_1 = C_2 = 0$.  Differentiating the structure
equations for $d\w^1$ and $d\w^2$
and reducing modulo $\theta_1, \w^1$, and $\w^2$ yields
\begin{align*}
0 = d(d\w^1) & \equiv (B_1 \theta_2\, - E_1\, \zeta)  \& \wt^3 \& \wt^4
\mod{\theta_1, \w^1, \w^2}\\
0 = d(d\w^2) & \equiv (B_2 \theta_2\, - E_2\, \zeta)  \& \wt^3 \& \wt^4
\mod{\theta_1, \w^1, \w^2};
\end{align*}
therefore $B_1 = B_2 = E_1 = E_2 = 0$.  A similar argument shows that
if $[C_3\  C_4]$
vanishes, then so do $[B_3\  B_4]$ and $[E_3\  E_4]$.
\end{pf}

In all the examples considered below, both vectors $[C_1 \ C_2]$ and
$[C_3 \ C_4]$
are nonzero.

The following two propositions are proved in \cite{C01b}.
\begin{prop}
If both vectors $[C_1 \ C_2], \ [C_3 \ C_4]$ vanish, then the \BB/
transformation
is locally equivalent to a transformation between solutions of the
wave equation $z_{xy} = 0$.
\end{prop}

\begin{prop}
If both of the vectors $[C_1 \ C_2], \ [C_3 \ C_4]$ are nonzero,
then the vectors $[B_1\ B_2], \ [B_3 \ B_4]$ are either both zero or
both nonzero.
\end{prop}

If the vectors $[B_1\  B_2]$ and $[B_3\  B_4]$ both
vanish, then for each $\lambda$ the \BB/ transformation is {\em holonomic}
in the sense described in \cite{C01b}.

\begin{prop}
If both of the vectors $[C_1 \ C_2], \ [C_3 \ C_4]$ are nonzero,
then the vectors $[D_1\  D_2], \ [D_3\  D_4]$
are either both zero or both nonzero.  Moreover, if these vectors both vanish,
then so do the vectors $[E_1 \ E_2], \ [E_3 \ E_4]$ and the functions
$F_1, \, F_2$.
\end{prop}
\begin{pf}
Suppose that $D_1 = D_2 = 0$. Differentiating the structure equation
for $d\theta_2$
and reducing modulo $\theta_1, \theta_2$, either $\w^1$ or $\w^2$,
and either $\w^3$ or $\w^4$ yields
\begin{alignat*}{2}
0 = d(d\theta_2) & \equiv &-C_1 D_4\, \zeta \& \w^2 \& \w^4
&\mod{\theta_1, \theta_2,\, \w^1, \w^3}\\
0 = d(d\theta_2) & \equiv & C_2 D_4\, \zeta \& \w^1 \& \w^4
&\mod{\theta_1, \theta_2,\, \w^2, \w^3}\\
0 = d(d\theta_2) & \equiv & -C_1 D_3\, \zeta \& \w^2 \& \w^3
&\mod{\theta_1, \theta_2,\, \w^1, \w^4}\\
0 = d(d\theta_2) & \equiv & C_2 D_3\, \zeta \& \w^1 \& \w^3
&\mod{\theta_1, \theta_2,\, \w^2, \w^4}.
\end{alignat*}
Since $C_1$ and $C_2$ are not both zero, we must have $D_3 = D_4 = 0$.
A similar argument shows the converse.

Now suppose that $D_1 = D_2 = D_3 = D_4 = 0$.
Differentiating the structure equation for $d\theta_1$
yields
\[0 = d(d\theta_1) \equiv
\zeta \& (E_3\, \theta_2 \& \w^4 - E_4\, \theta_2 \& \w^3 - F_1\, \w^3 \& \w^4)
\mod{\theta_1, \w^1, \w^2}; \]
therefore, $E_3 = E_4 = F_1 = 0$.  Similarly,
differentiating the structure equation for $d\theta_2$
yields
\[0 = d(d\theta_2) \equiv
\zeta \& (E_1\, \theta_1 \& \w^2 - E_2\, \theta_1 \& \w^1 - F_2\, \w^1 \& \w^2)
\mod{\theta_2, \w^3, \w^4}; \]
therefore, $E_1 = E_2 = F_2 = 0$.
\end{pf}

  From this proposition we see that
if either of the vectors $[D_1\  D_2], \
[D_3\  D_4]$ vanishes (and the corresponding $C$-vector is nonzero),
then the structure equations for $d\theta_i,\, d\w^i$
contain no terms involving $\zeta$.
In that case, the \BB/ transformation depends
trivially on the parameter $\lambda$.
Conversely, if these vectors are nonzero, then the \BB/ transformation
has nontrivial dependence upon $\lambda$ and may be considered
genuinely parametric.

\section{The sine-Gordon equation}\label{SGsec}
In this and subsequent sections, we will drop the usage of bars to
distinguish between objects on $\scP$ and the corresponding objects on $\G$.
(We will be working mostly on $\scP$ anyway.)  Instead, bars will be
used to distinguish
between similar coordinates (or similar coframes) on $\scM_1$ and $\scM_2$.

\subsection*{\BB/ Transformations}
As in the example given in the introduction,
let $\scM_1 = \scM_2 = \bb{R}^5$ with coordinates $(x,y,z,p,q)$ on $\scM_1$
and $(\bar{x},\bar{y},\bar{z},\bar{p},\bar{q})$
on $\scM_2$.  Let $\calI_1$ and $\calI_2$ be copies of the sine-Gordon system
generated by the forms
\begin{gather*}
    \theta_1 = dz - p\,dx - q\, dy \notag \\
    \Theta_1 = -dp \& dx - dq \& dy \\
    \Psi_1 = (dp - \tfrac{1}{2} \sin(2z)\, dy) \& dx \equiv -(dq - \tfrac{1}{2}
\sin(2z)\, dx) \& dy \mod \Theta_1 \notag
\end{gather*}
on $\scM_1$, and similar forms $\theta_2,\Theta_2,\Psi_2$ in terms of
the barred variables on $\scM_2$.
The \BB/ transformation $\scP \subset \scM_1 \times \scM_2 \times
\bb{R}$ is defined by the four equations
\begin{gather}
\bar{x} = x, \qquad \bar{y} = y \notag \\
p - \bar{p} = \lam\, \sin(z + \zbar)  \\
q + \bar{q} = \frac{1}{\lam}\, \sin(z - \zbar) \notag.
\label{SGPdef}\end{gather}

Let $\calK$ denote the span of the 1-forms $\theta_1$, $\theta_2$ and
$d\lambda$.  On $\scP$, we have
\begin{alignat*}{2}d\pbar - \tfrac{1}{2}\sin (2\zbar) \,dy &\equiv dp
- \tfrac{1}{2} \sin (2z)\, dy & \mod &{dx, \calK} \\
d\qbar - \tfrac{1}{2}\sin (2\zbar)\, dx &\equiv -(dq -\tfrac{1}{2}
\sin (2z) dx) & \mod &{dy, \calK}.
\end{alignat*}
  From these, it follows that $\Theta_1 \equiv -2\Psi_2$ and
$\Theta_2 \equiv -2\Psi_1$ modulo $\calK$.

\subsection*{$G$-structure Invariants}
We will use $x,y,z,\zbar, \pbar, q,$ and $\lambda$ as coordinates
on $\scP$; then $p, \qbar$ are given by the
\BB/ transformation equations
\begin{align}
p & = \pbar + \lambda\, \sin(z + \zbar)\label{SGBtrans} \\
\qbar & = -q + \frac{1}{\lam}\, \sin(z - \zbar).\notag
\end{align}
The coframe
$$
\begin{aligned}\zeta &=\lambda^{-1}d\lambda,\\
\theta_1 &= dz - p\,dx - q\, dy,\\
\theta_2 &= d\zbar-\pbar\, dx - \qbar\, dy,
\end{aligned}\qquad
\begin{aligned}
\w^1 &=dx,\\ \w^2 &= d\pbar- \tfrac{1}{2} \sin (2\zbar) \, dy,\\
\w^3 &=dy,\\ \w^4 &=dq - \tfrac{1}{2} \sin (2z)\,dx
\end{aligned}
$$
satisfies
\bel{dtzero}{
\begin{aligned}
d\theta_1 &\equiv A_1 \w^1 \& \w^2 + \w^3 \& \w^4\\
d\theta_2 &\equiv \w^1 \& \w^2 + A_2 \w^3 \& \w^4
\end{aligned}\qquad \mod \calK
}
with $A_1=1,\ A_2=-1$.  To get a section of $\scG$ we can modify
$\w^2$ and $\w^4$ to
\begin{align*}
\w^2 &= d\pbar- \tfrac{1}{2} \sin (2\zbar)\,
dy+\lambda\, (\sin(z + \zbar)\,\zeta + \cos(z + \zbar)\,\theta_2)\\
\w^4 &= dq - \tfrac{1}{2} \sin (2z)\,dx
+\lambda^{-1} (\sin(z - \zbar)\,\zeta - \cos(z - \zbar)\,\theta_1).
\end{align*}
Using this modified coframe, we compute that
\begin{alignat*}{3}
\begin{bmatrix} B_1 & C_1 & E_1\\
   B_2 & C_2 & E_2 \end{bmatrix} &=
   -\lam \begin{bmatrix} 0 & 0 & 0\\
    \sin (z + \zbar) & \cos (z + \zbar) & \cos (z + \zbar)
\end{bmatrix} \qquad &
\begin{bmatrix} D_1 \\ D_2 \end{bmatrix} &=
\lam\begin{bmatrix}\sin(z-\zbar) \\ 0 \end{bmatrix} \\[0.1in]
\begin{bmatrix} B_3 & C_3 & E_3 \\
B_4 & C_4 & E_4 \end{bmatrix} &=
\dfrac{1}\lam \begin{bmatrix} 0 & 0 & 0\\
\sin (z - \zbar) & -\cos (z - \zbar) & \cos (z - \zbar)\end{bmatrix}\qquad &
\begin{bmatrix} D_3 \\ D_4 \end{bmatrix} &=
\lam^{-1}\begin{bmatrix} \sin(z-\zbar)\\0 \end{bmatrix}.
\end{alignat*}
As well, $F_1=F_2=0$.  We observe that
\begin{enumerate}
\item{The vector triples $\{[ B_1 \  B_2],\ [C_1\ C_2],\ [ E_1 \  E_2]\}$ and
$\{[ B_3 \ B_4],\ [C_3\ C_4],\ [ E_3 \  E_4]\}$ each span a one-dimensional
space.}\label{lindepcond}
\item{Each $D$-vector is perpendicular to the corresponding
$C$-vector.}\label{perpcond}
\item{For both $\calI_1$ and $\calI_2$, each characteristic system
contains a rank one integrable subsystem, namely, those spanned
by $\{dx\}$ and $\{dy\}$.}\label{intsubsystemcond}
\end{enumerate}
In section \ref{concon} we will show that condition
(\ref{intsubsystemcond}) follows from conditions (\ref{lindepcond})
and (\ref{perpcond}).

\subsection*{Symmetries}
The sine-Gordon equation has a three-dimensional group of symmetries,
generated by translations in $x$ and $y$ and the
{\em Lie transformation} $L_\mu$, defined by
$$x \mapsto \mu^{-1} x,\quad y \mapsto \mu y, \quad z \mapsto
z,\qquad \mu \ne 0.$$
This is a symmetry of the sine-Gordon equation \eqref{SGE}, as well
as any ``f-Gordon''
equation $z_{xy}=f(z)$.  In other words, $L_\mu$ takes the graph of
one solution $z$
to that of a new solution $z'=L_\mu \cdot z$, defined by
$$z'(x,y)=z(\mu x, \mu^{-1} y).$$
Now, suppose $z,\zbar$ are two solutions that satisfy the \BT/
equations \eqref{SGEback}
for $\lambda=\lambda_0$.  Let the Lie transformation act
simultaneously on both, obtaining
new solutions $z',\zbar'$ respectively.  Then $z',\zbar'$ satisfy
\eqref{SGEback} for
$\lambda=\mu\lambda_0$.  In other words, the parameter in
\eqref{SGEback} is generated by
lifting the Lie transformation.

This observation is not new; in fact, Rogers and Shadwick \cite{rs}
show that, assuming that
the \BT/ does not change the independent variables, and that the
symmetries of each differential
equation are {\em point transformations} which cover the same
transformation of the
independent variables (e.g., $x \mapsto \mu^{-1} x$, $\  y \mapsto
\mu y$), then this
procedure always produces a one-parameter family of \BB/
transformations from a single
\BT/.  While it would be interesting to generalize this result to a
less restrictive setting,
at present we are concerned with how simultaneous symmetries in each
of the \MA/ equations
manifest themselves as diffeomorphisms of the 7-dimensional manifold
$\scP$.  Furthermore,
we will focus on infinitesimal symmetries, i.e., the vector fields
which generate symmetry
transformations.

A vector field $\vv$ is an (infinitesimal) symmetry of an EDS $\I$ if
$\Lie_\vv \psi\in \I$ for any form $\psi\in \I$.
(Here $\Lie$ denotes the Lie derivative.)
Such vector fields form a Lie algebra $\fs_\I$ under the usual bracket.  A
special subalgebra are the {\it Cauchy characteristic} vector fields,
which are simply vector fields annihilated by the 1-forms in the
Cartan system of $\I$.  This subalgebra $\fc_\I$ is either trivial
or infinite-dimensional, since a Cauchy characteristic vector field
may be multiplied by a smooth function.

It is straightforward to check that
$$\vv=x \dib{x}-y\dib{y} -p\dib{p}+q\dib{q}$$
is the infinitesimal symmetry of the \MA/ system \eqref{SGEMA}
generating the Lie transformation.
Let $\vvbar$ denote the corresponding vector field on $\scM_2$.
\begin{prop}\label{nocommonlift}
In terms of coordinates $x,y,p,q,z,\zbar$ and $t=\log \lambda$ on $\scP$,
$$\ww=x \dib{x}-y\dib{y} -p\dib{p}+q\dib{q}-\dib{t}$$
is the unique vector field on $\scP$ which covers $\vv$ and $\vvbar$, i.e.,
$\pi_{1*}\ww=\vv$ and $\pi_{2*}\ww=\vvbar$.
\end{prop}
\begin{proof} Suppose that
$$\ww =x \dib{x}-y\dib{y} -p\dib{p}+q\dib{q}+A\dib{\zbar}+B\dib{t}$$
covers $\vv,\vvbar$.
Since $\vvbar\intprod d\zbar=0$, we have $A=0$.  Since $\vvbar\intprod
d\pbar = -\pbar$,
we must have
\begin{align*}
-\pbar & = \ww \intprod d\pbar \\
& = \ww \intprod [dp - \lambda\cos(z+\zbar)(dz+d\zbar)
-\lambda\sin(z+\zbar)\, dt]\\
& = -p - \lambda \sin(z + \zbar)\, B.
\end{align*}
Then \eqref{SGBtrans} gives $B=-1$.
\end{proof}

Note that Prop. \ref{nocommonlift}
   implies that there is no vector field tangent to the 6-dimensional
level sets $\B_\lambda$ of
$\lambda$ which covers $\vv$ and $\vvbar$.  However, as has been
observed by Igonin and
Krasil'shchik \cite{igback}, there is a lift of $\vv$ to $\B_1$ which
generates the one-parameter
family of \BB/ transformations.  We will now make this construction explicit.

Let $\vvhat$ be the lift of $\vv$ to $\B=\B_1$ such that $\vvhat
\intprod d\zbar=0$.  Then $\vvhat$
is not a symmetry of the form $\theta_2$, since
$\Lie_\vvhat\theta_2= \sin(z+\zbar)dx + \sin(z-\zbar)dy.$
Thus, flow by $\vvhat$ induces a 1-parameter group of diffeomorphisms
$\varphi_t$ of $\B$
which preserve $\theta_1$ up to
multiple (i.e., $\varphi_t^*\theta_1\equiv 0 \mod \theta_1$)
but don't preserve $\theta_2$.  These are
$$\varphi_t(x,y,p,q,z,\zbar)=(e^t x,e^{-t}y,e^{-t}p,e^t q,z,\zbar),$$

Now let $\varphi:\B\times \R\to \B$ be defined by
$\varphi(x,y,p,q,z,\zbar,t)=\varphi_t(x,y,p,q,z,\zbar)$.
One computes that
$$\varphi^* \theta_2=d\zbar - (p-\lambda \sin(z+\zbar))dx -
(-q+\lambda^{-1}\sin(z+\zbar))dy
\mod dt,$$
where we have set $\lambda=e^t$.  Thus, by splitting off the $dt$
parts of these forms and
re-labeling $\varphi^*\theta_2$ as $\theta_2$ on $\scP=\B\times \R$,
we recover the parametric
\BT/ defined earlier.

Note that it is not clear how to choose the lift $\vvhat$ so as to
make this construction
work, nor is this question addressed in \cite{igback}.
However, some light may be shed in \S\ref{Zsec},
where we will see how this trick applies in another example.

\subsection*{Conservation laws}
The space of {\em conservation laws} for an EDS $\calI$ whose
integral manifolds have dimension $2$ may be naturally identified
with the set of closed $2$-forms in $\calI$, modulo the exterior
derivatives of $1$-forms contained in $\calI$, i.e., with the quotient space
$$\frac{\{\Phi \in \calI^2\,\vert\, d\Phi = 0\}}{\{d\phi\,\vert\,
\phi \in \calI^1\}}.$$
Any conservation law for the Monge-Amp\`ere system $\calI_i =
\{\theta_i,\, \Theta_i,\, \Psi_i\}$ has a unique
representative of the form
\begin{equation}
   \Phi = Q\, \Psi_i + \theta_i \& \gamma \label{MACLform}
\end{equation}
for some function $Q$ and $1$-form $\gamma$, while any conservation
law for the system $\calJ = \{\theta_1,\, \theta_2,\,
\Theta_1,$ $\Theta_2\}$ on $\scB$ has a unique representative of the
form
\begin{equation}
   \Phi = \theta_1 \& \gamma_1 + \theta_2 \& \gamma_2  \label{BTCLform}
\end{equation}
for some $1$-forms $\gamma_1,\, \gamma_2$ on $\scB$.  (It is not
difficult to show that any conservation law for the system
$\calJ$ on any slice $\scB_{\lambda}$ has a unique lift to a
conservation law for the system $\calK = \{\theta_1,\,
\theta_2,\, \zeta\}$ on $\scP$, and that conversely, any
conservation law for $\calK$ restricts to a conservation law
for $\calJ$ on each slice $\scB_{\lambda}$. Thus, for ease of computation
we will work on $\scB$ rather than on $\scP$ when
computing conservation laws.)  Since $\pi_i^*\calI_i \subset \calJ$,
the pullback of any conservation law for either of the Monge-Amp\`ere
systems $\calI_i$ to $\scB$ is also a conservation law for
$\calJ$, though it typically must be modified by the exterior
derivative of a 1-form in $\calJ$ in order to appear in the
form \eqref{BTCLform}.

The condition that $\Phi$ be a closed form in the ideal typically
leads to an overdetermined system of PDE's which must be satisfied by
the coefficients of $\Phi$.  Here we will sketch the computation of
the spaces of conservation laws in the sine-Gordon example.

First we compute the conservation
laws for the ideal $\calI_1$.  (This
will, of course, be isomorphic to
the space of conservation laws for
$\calI_2$.)  Suppose that
\[ \Phi = Q\,[(dp - \tfrac{1}{2} (\sin z)\, dy) \& dx - (dq -
\tfrac{1}{2} (\sin z)\, dx) \& dy] + \theta_1 \& \gamma \]
is a closed form in $\calI_1$.  Computing
$d\Phi \equiv 0 \mod{\theta_1}$
shows that
\[ \gamma = -Q_p\, dp + Q_q\, dq - (Q_x +p\, Q_z + \sin(2z)\, Q_q)\,
dx + (Q_y +q\, Q_z + \sin(2z)\, Q_p)\, dy. \]
Then the condition $d\Phi = 0$ gives a system of second-order PDE's
for the function $Q$ whose solutions are
\begin{equation}
   Q = c_1\, (px - qy) + c_2\, p + c_3\, q \label{SGECLfct}
\end{equation}
for arbitrary constants $c_1, c_2, c_3$.  Thus the space of
conservation laws for $\calI_1$ (and hence for $\calI_2$ as well) is
three-dimensional.

\begin{rem} By applying the method of equivalence to hyperbolic \MA/ systems,
one can determine which systems are {\em variational}, i.e. are
contact-equivalent
to the Euler-Lagrange equation for a first-order Lagrangian $\int
L(x,y,u,p,q)\,dx\,dy$.
For those which are variational, N\"oether's theorem gives a
one-to-one correspondence
between conservation laws and symmetries for the \MA/ system.
The sine-Gordon equation is variational and has a three-dimensional
symmetry group, so
it is no surprise that we arrive at a three-dimensional space of
conservation laws.
\end{rem}

Computing the space of conservation laws for $\calJ$ is
considerably more involved.  We will use the coframing
described above (setting $\zeta=0$ because we are restricting to
$\scB$), i.e.,
\begin{align*}
\theta_1 &= dz - p\,dx - q\, dy\\
\theta_2 &= d\zbar-\pbar\, dx - \qbar\, dy \\
\w^1 &=dx\\
\w^2 &= d\pbar- \tfrac{1}{2} \sin (2\zbar)\,
dy+\lambda\, \cos(z + \zbar)\,\theta_2\\
\w^3 &=dy\\
\w^4 &= dq - \tfrac{1}{2} \sin (2z)\,dx
-\lambda^{-1}\, \cos(z - \zbar)\,\theta_1.
\end{align*}
Suppose that
\[ \Phi = \theta_1 \& (P_1\, \w^1 + P_2\, \w^2 + P_3\, \w^3 + P_4\,
\w^4) + \theta_2\& (Q_1\, \w^1 + Q_2\, \w^2 + Q_3\, \w^3 + Q_4\,
\w^4) + R\, \theta_1 \& \theta_2 \]
is a closed form in $\calJ$.  Computing
$$d\Phi \equiv 0 \mod{\theta_1, \theta_2}$$
shows that
$Q_1 = P_1$, $Q_2 = P_2$, $Q_3 = -P_3$, and $Q_4 = -P_4$.
Next, computing
\[ d\Phi \equiv 0 \mod{\w^1, \w^2, \theta_1 - \theta_2} \]
shows that
\[ R = \frac{1}{\lambda} \cos (z - \zbar) P_4 - \lambda \cos(z + \zbar) P_2. \]
Then the condition $d\Phi=0$ gives a system of 14 first-order PDE's
for the four functions $P_1, P_2, P_3, P_4$.  The compatibility
conditions for this system lead to five additional equations; the
resulting system is not involutive, so it must be prolonged.  Further
analysis shows that the general solution of this system is
\begin{align*}
P_1 & = c_1[\pbar + \tfrac{1}{2} \lam \sin (z + \zbar) - \tfrac{1}{2}
y\sin(2z) + \lam(\pbar x + q y)\cos(z + \zbar)] + c_2 \lambda\, \pbar
\cos (z + \zbar) \\
&  \qquad  + c_3 (\lambda\, q \cos (z + \zbar) - \tfrac{1}{2} \sin
(2z)) + \tfrac{1}{2} c_4 \lambda \sin (z + \zbar) \\
P_2 & = c_1 x + c_2 \\
P_3 & = c_1[q - \tfrac{1}{2}\lam^{-1}  \sin (z - \zbar) +
\tfrac{1}{2} x\sin(2\zbar) + \lam^{-1}(\pbar x - q y)\cos(z - \zbar)]
\\
& \qquad  + c_2 (\lambda^{-1} \pbar \cos (z - \zbar) + \tfrac{1}{2}
\sin (2\zbar)) - c_3 \lambda^{-1} q \cos (z - \zbar) + 
\tfrac{1}{2} c_4 \lambda^{-1} \sin (z - \zbar) \\
P_4 & = c_1 y + c_3
\end{align*}
for arbitrary constants $c_1, c_2, c_3, c_4$.  Thus the space of
conservation laws for $\calJ$ is four-dimensional.  (See \cite{EDS} for
details about Cartan-K\"ahler analysis of exterior differential systems.)

Consider the following question: given a conservation law for
$\calI_1$, is its pullback to $\calJ$ equivalent to the
pullback of a conservation law for $\calI_2$?  In other words, do the
pullbacks of the spaces of conservation laws for $\calI_1,\, \calI_2$
span the same three-dimensional subspace of the conservation laws for
$\calJ$?  Perhaps surprisingly, the answer in this case is
no: it can be shown that any conservation law in $\calI_1$
corresponding to a function $Q$ of the form \eqref{SGECLfct} with
$c_1 \neq 0$ is not equivalent to the pullback of a conservation law
for $\calI_2$, and vice-versa.
Thus, the four-dimensional space of conservation laws for
$\calJ$ is the union of the pullbacks of the
three-dimensional spaces of conservation laws for $\calI_1$ and
$\calI_2$.

\section{Goursat's Equation}\label{Zsec}
In an 1899 memoir, Goursat \cite{G1899} classified those equations of
the form $u_{xy}=f(x,y,u,u_x,u_y)$ which are integrable by the method
of Darboux (without prolongation).  The equation
\bel{zvy}{w_{xy}+2\dfrac{\sqrt{w_x}\sqrt{w_y}}{x+y}=0}
occurs at the head of Goursat's list, and we will refer to it here as
`Goursat I' or simply as Goursat's equation.
(We will restrict our attention to solutions $w(x,y)$ for which
the arguments of the square roots in \eqref{zvy} are positive.)

\subsection*{B\"acklund transformations and Symmetries}
A \BT/ linking solutions of \eqref{zvy} with those of the
wave equation $u_{xy}=0$ was given by Zvyagin \cite{zvyagin}, and takes the
form
\bel{ztrans}{u_x=\dfrac{u}{x+y}+\sqrt{2w_x},\qquad
u_y=\dfrac{u}{x+y}+\sqrt{2w_y}.}
As in the sine-Gordon example we may use the symmetries of
\eqref{zvy} to construct a
parametric version of \eqref{ztrans}.
Let $\scM_1=J^1(\R^2,\R)$ with coordinates $x,y,w,p,q$, and
$\scM_2=J^1(\R^2,\R)$ with
coordinates $\xbar,\ybar,u,\pbar,\qbar$.  We regard \eqref{ztrans}, together
with $x=\xbar$ and $y=\ybar$, as defining $\B^6 \subset \scM_1\times
\scM_2$.  We pull back the
contact forms to $\B$ to give the usual rank two Pfaffian system generated by
$$\theta_1=dw-p\,dx-q\,dy,\qquad \theta_2=du-\pbar dx - \qbar dy,$$
where
$$\pbar=\dfrac{u}{x+y}+\sqrt{2p},\qquad \qbar=\dfrac{u}{x+y}+\sqrt{2q}.$$

In an earlier discussion of Darboux integrability (\cite{Gvol2}, p.
196), Goursat
cited \eqref{zvy} as one of the few known
Darboux-integrable equations admitting a finite-dimensional
symmetry group.  In fact, its symmetries take the form
$$w\mapsto A w+B,\qquad x\mapsto\dfrac{ax+b}{cx+d}, \qquad
y\mapsto-\dfrac{ay-b}{cy-d}$$
where $A\ne 0$ and we may take
$\left(\begin{smallmatrix}a & b\\c & d\end{smallmatrix}\right)\in SL(2,\R)$.
In the case of linear substitution for $w$, or for $x$ and $y$, the
vector field
on $\scM_1$ which generates the symmetry has a lift to $\B$ which also
covers a symmetry
vector field on $\scM_2$.  However, consider the vector field
$$\vv=x^2\dib{x}-y^2\dib{y}-2xp\dib{p}+2yq\dib{q},$$
on $\scM_1$, which generates the following 1-parameter group of
symmetries for \eqref{zvy}:
\bel{zvact}{x\mapsto \dfrac{x}{1-tx},\quad y\mapsto\dfrac{y}{1+ty}, \quad
p\mapsto (1-tx)^2p,\quad q\mapsto (1+ty)^2 q,\quad w\mapsto w.}

\begin{prop} There is no lift of $\vv$ to $\B$ which covers a
symmetry of the wave equation
on $\scM_2$.
\end{prop}
\begin{proof} Suppose
$$\vvhat=x^2\dib{x}-y^2\dib{y}-2xp\dib{p}+2yq\dib{q}+f\dib{u}$$
is such a lift, where $f$ is an unknown function of $x,y,w,p,q,u$.
The condition
$\Lie_\vvhat \theta_2 \equiv 0 \mod \theta_2$ implies that
$$df=g\,\theta_2 + \left(\dfrac{f}{x+y}+u+x\sqrt{2p}\right)dx +
\left(\dfrac{f}{x+y}-u-x\sqrt{2q}\right)dy,$$
where $g$ is necessarily a function of $x,y$ and $u$ only.  Further
differentiation
gives
$$0=d^2f\equiv \dfrac{x-g}{\sqrt{2p}} dp \& dx -\dfrac{y+g}{\sqrt{2q}} dq \& dy
\mod du \& dx, du \& dy, dx \& dy,$$
implying that $g=x$ and $g=-y$, a contradiction.
\end{proof}

Now let $\varphi_t$ be the 1-parameter family of diffeomorphisms of
$\B$ defined by
\eqref{zvact} along with $u\mapsto u$.  Let $\tilde\scP\subset \B
\times \R$ be the open set
on which $\varphi_t$ is defined, and let $\varphi:\tilde\scP \to \B$
be defined by
$\varphi(x,y,p,q,u,w,t)=\varphi_t(x,y,p,q,u,w)$.  Then one computes that
$$\varphi^*\theta_2 \equiv du -
\left(\dfrac{u(1+ty)}{x+y}+\sqrt{2p}\right)\dfrac{dx}{1-tx}
-\left(\dfrac{u(1-tx)}{x+y}+\sqrt{2q}\right)\dfrac{dy}{1+ty}\mod dt.$$
This form is the same as the pullback of $\theta_2$ to the
submanifold $\scP \subset \scM_1 \times \scM_2
\times \R$ defined by
\bel{zbtt}{(1-tx)\pbar =\dfrac{u(1+ty)}{x+y}+\sqrt{2p},\qquad
(1+ty)\qbar =\dfrac{u(1-tx)}{x+y}+\sqrt{2q}.}
It is straightforward to verify that these equations define a \BT/
for any $t$ for which
the maps \eqref{zvact} are defined.

\subsection*{$G$-Structure Invariants}
The coframe
\begin{equation}\label{zframe}
\begin{aligned}\zeta&=\left((1-tx)(1+ty)\right)^{-1}dt\\
\theta_1 &=\frac{1}{\sqrt{2q}}\,(dw-p\,dx-q\,dy)\\
\theta_2 &=(1-tx)(du-\pbar\,dx-\qbar\,dy)
\end{aligned}\qquad
\begin{aligned}
\w^1 & = dx\\
\w^2 &=\dfrac{dp}{\sqrt{2p}}+\dfrac{\sqrt{2q}}{x+y}dy\\
\w^3 & = dy\\
\w^4 &=\dfrac{dq}{\sqrt{2q}}+\dfrac{\sqrt{2p}}{x+y}dx
-\left(\dfrac{u+y\sqrt{2q}}{1+ty}\right) dt
\end{aligned}\end{equation}
gives a section of $\G$ with $A_1=\sqrt{p/q}$ and $A_2=(1-tx)/(1+ty)$.
With respect to this coframe, we compute that
\begin{alignat*}{3}
\begin{bmatrix} B_1 & C_1 & E_1\\
   B_2 & C_2 & E_2 \end{bmatrix} &=
\begin{bmatrix} 0 & 0 & 0\\
                  0 & -(x+y)^{-1}& 0 \end{bmatrix} \qquad &
\begin{bmatrix} D_1 \\ D_2 \end{bmatrix} &=
\begin{bmatrix}-(1+ty)(u+x\sqrt{2p}) \\ 0 \end{bmatrix} \\[0.1in]
\begin{bmatrix} B_3 & C_3 & E_3 \\
B_4 & C_4 & E_4 \end{bmatrix} &=
\begin{bmatrix} 0 & 0 & 0\\
0 & -(x+y)^{-1} &-1\end{bmatrix}\qquad &
\begin{bmatrix} D_3 \\ D_4 \end{bmatrix} &=
\begin{bmatrix} -(1-tx)(u+y\sqrt{2q})\\0 \end{bmatrix}.
\end{alignat*}
In this case, neither $F_1$ nor $F_2$ is zero, nor is their ratio a
constant.  We observe that
\begin{enumerate}
\item{Both vectors $[B_1\ B_2]$ and $[B_3\ B_4]$ are zero, indicating
that, for each value of
$t$, this \BT/ is holonomic (see \cite{C01b}, Thm. 6.2).}
\item{Both sets of $B$-, $C$- and $E$-vectors are linearly dependent,
which is consistent
with the presence of one-dimensional integrable subsystems in the
characteristic systems.}
\item{As in the sine-Gordon example, both $D$-vectors are
perpendicular to the corresponding
$E$-vectors.}
\end{enumerate}

\begin{rem} While the transformation \eqref{ztrans} is implicit in
\cite{zvyagin},
the transformation given there is actually
\bel{zothertrans}{z_x=\left(\sqrt{w_x}+\sqrt{\dfrac{z-w}{x+y}}\right)^2,\qquad
   z_y=\left(\sqrt{w_y}+\sqrt{\dfrac{z-w}{x+y}}\right)^2,}
where $z(x,y)$ satisfies the wave equation $z_{xy}=0$ if and only if
$w(x,y)$ satisfies \eqref{zvy}.  Moreover, $z$ is related to $u$
via the change of variable $z=w+\frac12 u^2/(x+y)$, implying that
\bel{zutrans}{z_x=\frac12(u_x)^2, \qquad z_y=\frac12(u_y)^2.}
In other words, the two wave equation solutions $z$ and $u$
are related by the simple \BT/ \eqref{zutrans}.  Such transformations
can be characterized by the vanishing of all $B$- and $C$-vectors
(see \cite{C01b}, Thm. 3.1).

An interesting observation is that the \BT/ \eqref{zothertrans} is
clearly the composition
of the transformations \eqref{ztrans} and \eqref{zutrans}, each of
which exhibit
simple values for the $B$- and $C$-vector invariants (e.g. all the $B$-vectors
are zero).  However, one can compute that for \eqref{zothertrans}
the $B$-vectors are nonzero and linearly dependent on the $C$-vectors.
\end{rem}

\subsection*{Conservation laws}
We will compute the spaces of conservation laws for the wave equation,
for Goursat's equation \eqref{zvy},
and for the ideal $\calJ$ on $\scB_t$ for any fixed value of $t$.

First we compute the conservation laws for the system
\[ \calI_2 = \{du - \pbar\, dx - \qbar \, dy, \ d\pbar \& dx, \
d\qbar \& dy\}, \]
which represents the wave equation $u_{xy}=0$.  Suppose that
\[ \Phi = Q\,[d\pbar \& dx - d\qbar \& dy] + \theta_2 \& \gamma \]
is a closed form in $\calI_2$.  Computing
$d\Phi \equiv 0 \mod{\theta_2}$
shows that
\[ \gamma = -Q_{\pbar}\, d\pbar + Q_{\qbar}\, d\qbar - (Q_x +\pbar\, Q_z)\,
dx + (Q_y +\qbar\, Q_z)\, dy. \]
Then the condition $d\Phi = 0$ gives a system of second-order PDE's
for the function $Q$ whose solutions are
\begin{equation}
   Q = f(x,\pbar) + g(y,\qbar) \label{waveCLfct}
\end{equation}
for arbitrary functions $f,g$.  Thus the space of
conservation laws for $\calI_2$ is infinite-dimensional and depends
(in the sense of Cartan-K\"ahler) on two arbitrary functions of two variables.

Next we compute the conservation laws for the system
\[ \calI_1 = \{dw - p\, dx - q \, dy, \
\left(dp + \frac{2\sqrt{pq}}{x+y}\, dy\right) \& dx, \
\left(dq + \frac{2\sqrt{pq}}{x+y}\, dx\right) \& dy\}, \]
which represents Goursat's equation \eqref{zvy}.
This computation is most easily carried out using the coframe
\begin{align*}
& \qquad \ \  \theta_1 = dw - p\, dx - q \, dy, \\[0.1in]
& \begin{aligned} \w^1 & = dx, \\
\w^2 & = dp + \frac{2\sqrt{pq}}{x+y}\, dy
\end{aligned}\qquad
\begin{aligned}
\w^3 & = dy, \\
\w^4 & = dq + \frac{2\sqrt{pq}}{x+y}\, dx.
\end{aligned}
\end{align*}
Suppose that
$ \Phi = Q\,[\w^1 \& \w^2 - \w^3 \& \w^4] + \theta_1 \& \gamma$
is a closed form in $\calI$.  Setting
\[ dQ = Q_0\, \theta_1 + Q_1\, \w^1 + Q_2\, \w^2 + Q_3\, \w^3 + Q_4\, \w^4 \]
and computing
$d\Phi \equiv 0 \mod{\theta_1}$
shows that
\[ \gamma = \left(Q_1 - \frac{2\sqrt{p}}{(x+y)\sqrt{q}}\,Q\right)\, \w^1
+ Q_2\, \w^2 - \left(Q_3 -
\frac{2\sqrt{q}}{(x+y)\sqrt{p}}\,Q\right)\, \w^3 - Q_4\, \w^4. \]
Then the condition $d\Phi = 0$ gives the first-order PDE
\[ Q_0 = \frac{2p\,Q_2 + 2q\,Q_4 + Q}{2(x+y)\sqrt{pq}} \]
and four additional second-order PDE's for $Q$.
(Some care must be taken here, as mixed partial derivatives in terms
of this coframing do not commute.) The compatibility
conditions for this system lead to three additional
equations; the resulting system is not involutive, so it
must be prolonged.  The prolonged system is involutive
with last nonvanishing Cartan character $s_1=2$, so the
space of solutions (and hence the space of conservation
laws for \eqref{zvy}) depends on two arbitrary functions
of one variable.  Note that \eqref{zvy} is not
variational, so there is no contradiction between the
finite-dimensional symmetry group and the
infinite-dimensional space of conservation laws.

Finally, we compute the space of conservation laws for the system
\[ \calJ = \{\theta_1,\, \theta_2,\, \Theta_1,\, \Theta_2\} \]
on $\scB_t$ for some fixed $t \in \bb{R}$.  We will use the coframe
$(\theta_1,\theta_2,\w^1,\w^2,\w^3,\w^4)$ defined by \eqref{zframe},
restricted to $\scB$.

Suppose that
\[ \Phi = \theta_1 \& (P_1\, \w^1 + P_2\, \w^2 + P_3\, \w^3 + P_4\,
\w^4) + \theta_2\& (Q_1\, \w^1 + Q_2\, \w^2 + Q_3\, \w^3 + Q_4\,
\w^4) + R\, \theta_1 \& \theta_2 \]
is a closed form in $\calJ$.  Computing
$d\Phi \equiv 0 \mod{\theta_1, \theta_2}$
shows that
\begin{alignat*}{2}
(1 + ty)\, P_1 - (1 - tx)\, Q_1 &= 0, & \qquad
\sqrt{p}\, P_3 + \sqrt{q}\, Q_3 &= 0, \\
(1 + ty)\, P_2 - (1 - tx)\, Q_2 &= 0, & \qquad
\sqrt{p}\, P_4 + \sqrt{q}\, Q_4 & = 0 .
\end{alignat*}
Thus we may write
\begin{multline*}
\Phi = ((1 - tx)\, \theta_1+ (1 + ty)\, \theta_2) \& (S_1\, \w^1 +
S_2\, \w^2)\\
  +(\sqrt{p}\, \theta_1 - \sqrt{q}\, \theta_2) \& (S_3\, \w^3 + S_4\, \w^4) +
  R\, \theta_1 \& \theta_2
\end{multline*}
for some functions $S_1, S_2, S_3, S_4$.  Computing
$d\Phi \equiv 0 \mod{\w^1, \w^2, \sqrt{p}\, \theta_1 - \sqrt{q}\, \theta_2}$
shows that
\[ R = \frac{(1 + ty)\, S_2 + \sqrt{q}\, S_4}{x + y}. \]
Then the condition $d\Phi=0$ gives a system of 14 first-order PDE's
for the four functions $S_1, S_2, S_3, S_4$.  The compatibility
conditions for this system lead to three additional equations; the
resulting system is not involutive, so it must be prolonged.
The prolonged system is involutive with last nonvanishing Cartan
character $s_2=2$,
so the space of solutions (and hence the space of conservation laws
for $\calJ$)
depends on two arbitrary functions of two variables.

Clearly, not every conservation law in $\calJ$ is the pullback of a
conservation law
for \eqref{zvy}, since the space of conservation laws for $\calJ$ is
strictly larger than that for \eqref{zvy}.  It may be that every conservation
law in $\calJ$ is the pullback of a conservation law for the wave equation,
but because the Cartan-K\"ahler analysis does not give explicit solutions,
this hypothesis is difficult to confirm or disprove.

\section{Pseudospherical surfaces in Euclidean space}\label{pseudosec}
In this section we study pseudospherical surfaces as arising from integrals
of a \MA/ system.  We present the classical \BB/ transformation in the
same context, as a double fibration fitting the description in \S\ref{intro}.
We calculate the torsion coefficients for the corresponding $G$-structure,
and note how these differ from those in the sine-Gordon case.
Since we are also interested in cases where a one-parameter family of
\BB/ transformations
may be generated by a symmetry of one of the \MA/ systems involved,
we also calculate the
symmetries of the system for pseudospherical surfaces.
\subsubsection*{Surface geometry via moving frames}
Let $\F$ be the bundle of Euclidean orthonormal frames on $\R^3$, on
which the basepoint
projection is $\xx$ and the vector-valued functions $e_1,e_2,e_3$
give the members of
the frame at $\xx \in \R^3$.  On $\F$ there are {\em canonical 1-forms}
   $\eta^i,\eta^i_j$ defined by
resolving the differentials of these vectors in the frame:
$$d\xx = e_i \eta^i, \qquad de_i = e_j \eta^j_i,$$
summation being understood on repeated indices.  (Because the $e_i$ are
orthonormal, the $\eta^i_j$ form
a skew-symmetric $3\by3$ matrix of 1-forms.)
By differentiating these defining equations we deduce
the {\em structure equations} for $\F$:
$$d\eta^i = -\eta^i_j \& \eta^j,\qquad d\eta^i_j = -\eta^i_k\& \eta^k_j.$$

For later use, we note that $\F$ can be identified with a matrix Lie
group on which
the $\eta$'s are left-invariant forms.  For, we may embed the vectors of
the frame in a $4\by 4$ matrix
$$g = \begin{pmatrix} 1 & 0 \\ \xx & A \end{pmatrix} \in  \R^3 \times SO(3),$$
where the columns of $A$ are $e_1,e_2,e_3$ in order.  Then
$$g^{-1}dg = \begin{pmatrix} 0 & 0 & 0 & 0 \\
\eta^1 & 0 & \eta^1_2 & \eta^1_3\\
\eta^2 & \eta^2_1 & 0 & \eta^2_3\\
\eta^3 & \eta^3_1 & \eta^3_2 &0 \end{pmatrix} \in \R^3 \oplus \fso(3).$$
The action of this group on $\F$ given by left-multiplication covers the
Euclidean motions in $\R^3$.

Given a regular surface $S\subset \R^3$, it is standard to construct
{\it first-order adapted} framings along $S$, which are lifts into $\F$
\begin{center}
\begin{picture}(30,35)(0,0)
\put(5,5){\makebox(0,0){$S\hookrightarrow$}}
\put(24,5){\makebox(0,0){$\R^3$}}
\put(22,29){\makebox(0,0){$\F$}}
\put(3,11){\vector(1,1){12.5}}
\put(22,24){\vector(0,-1){14}}
\put(27,19){\makebox(0,0){$\scriptstyle{\xx}$}}
\end{picture}
\end{center}
such that $e_1,e_2$ are tangent and $e_3$ is normal to the surface.
(Note that such
framings are not unique, since they may be modified by rotating $e_1,e_2$
within the plane they span.)  These lifts are surfaces $\Sigma\subset \F$
to which $\eta^3$ restricts to be zero, and along which $\eta^1,\eta^2$ are
linearly independent.  It follows from the structure equations that
\bel{heq}{\eta^3_i\restr_\Sigma = h_{i1} \eta^1+ h_{i2}\eta^2}
for a symmetric $2\by2$ matrix of functions $h_{ij}$, which are the components
of the second fundamental form of $S$ relative to the basis $e_1,e_2$
\cite{babyoneill}.
The Gauss curvature $K$ and mean curvature $H$ of the surface are respectively
equal to the determinant and one-half the trace of the matrix $(h_{ij})$.
Consequently, the restricted forms on $\Sigma$ satisfy
\begin{gather*}
\eta^3_1 \& \eta^3_2 = K \eta^1 \& \eta^2\\
\eta^3_1 \& \eta^2-\eta^3_2 \& \eta^1 = 2H \eta^1 \& \eta^2.
\end{gather*}

\subsubsection*{The $K=-1$ system}
A {\em pseudospherical surface} in $\R^3$ is one with $K=-1$ at each
point.  For
such surfaces, any first-order adapted framing is an integral surface
of the EDS
$\I$ generated by the 1-form $\eta^3$ and the 2-forms $d\eta^3$ and $$\Psi =
\eta^3_1 \& \eta^3_2 + \eta^1 \& \eta^2.$$ This would satisfy the Defn.
\ref{MAdefn} for \MA/ systems, except that the underlying manifold $\F$ is
6-dimensional.  However, the Cartan system of $\I$ is
$\{\eta^1,\eta^2,\eta^3,\eta^3_1,\eta^3_2\}$.  The integral curves of this
Frobenius system are the fibers of a projection $\F\to \scrU^5$, where $\scrU$ is the
unit tangent bundle of $\R^3$ and the projection is given by taking
$e_3$ as the
unit vector.  (Consequently, moving along the fibers corresponds to fixing the
basepoint of a frame in $\F$ and rotating $e_1$ and $e_2$.) By a standard
argument, the generating differential forms project to become
well-defined (up to
multiples) on $\scrU$, and so there is a well-defined EDS on $\scrU$ such that each
integral surface in $\F$ which is transverse to the fibers corresponds to a
unique integral surface in $\scrU$.  In this case, the EDS on $\scrU$ is a genuine \MA/
system which we will also call $\I$.  That this system is hyperbolic
follows from
the factorizations
$$\Psi+d\eta^3=(\eta^3_1+\eta^2) \& (\eta^3_2-\eta^1),\qquad
\Psi-d\eta^3 = (\eta^3_1-\eta^2) \& (\eta^3_2 + \eta^1).$$
In general, we may perform computations with
the generators of $\I$ up on $\F$, knowing that our conclusions will
hold on $\scrU$.

\subsubsection*{Relationship with sine-Gordon}
By breaking the Cauchy characteristic symmetry of the system for
pseudospherical surfaces, we can show how such surfaces may be constructed from
solutions of the sine-Gordon equation \eqref{SGE}.

On an open subset of the surface which is free of umbilic points,
we may choose a Darboux framing, i.e.
a first-order adapted framing which diagonalizes the second fundamental form.
For such a framing, we have
$$\begin{pmatrix}\eta^3_1 \\ \eta^3_2\end{pmatrix}=
\begin{pmatrix} \tan\phi & 0 \\ 0 & -\cot \phi \end{pmatrix}
\begin{pmatrix}\eta^1 \\ \eta^2\end{pmatrix},$$
where $\phi \in (0, \pi/2)$ is one-half of the angle between the
asymptotic lines.
We adjoin $\phi$ as a new variable, and enlarge our EDS to a rank
three Pfaffian system
$\{\eta^3,\eta^3_1-(\tan\phi)\, \eta^1, \eta^3_2+(\cot\phi)\,\eta^2\}$
on $\F\times \R$.  The vanishing of the exterior derivatives
of the last two 1-forms implies that the forms
$\eta^1/\cos\phi$ and $\eta^2/\sin\phi$ are closed.
Hence, on any integral surface there exist local
coordinates $t_1,t_2$ such that $\eta^1=(\cos\phi)\, dt_1$ and
$\eta^2=(\sin\phi)\, dt_2$.
It follows that $x=(t_1+t_2)/2$ and $y=(t_1-t_2)/2$ are arclength
coordinates along the asymptotic lines.

Substituting $d\phi=\phi_x\, dx + \phi_y\, dy$ into the 2-forms of
the EDS shows that
$\eta^1_2 = \phi_y\, dy - \phi_x\, dx$ along any integral surface.
It is then easy to see that the
structure equation $d\eta^1_2 = \eta^3_1 \& \eta^3_2$ implies that
$\phi$ satisfies the sine-Gordon equation
   $$\phi_{xy}=\tfrac{1}{2} \sin(2\phi).$$
Conversely, if $\phi(x,y)$ is a solution to the sine-Gordon equation,
we may produce a pseudospherical surface by integration.  For, if we
let $A=(e_1,e_2,e_3)$
as before, then the structure equations $de_i = e_j\eta^j_i$ imply that
\bel{Aeq}{\dfrac{\di A}{\di x}=
A \begin{pmatrix} 0 & -\phi_x & -\sin\phi \\ \phi_x & 0 & \cos \phi
\\ \sin\phi & -\cos\phi & 0\end{pmatrix}
\qquad
\dfrac{\di A}{\di y}=
A \begin{pmatrix} 0 & \phi_y & -\sin\phi \\ -\phi_y & 0 & -\cos \phi
\\ \sin\phi & \cos\phi & 0\end{pmatrix}.
}
(This is an overdetermined system for matrix $A(x,t)$,
and its integrability condition is the sine-Gordon equation for $\phi$.)
So, given a solution to sine-Gordon, we may obtain the framing by
solving linear systems of ODE,
and then solve for the surface
$\xx(x,y) \in \R^3$ by integrating the equations
\bel{xeq}{\dfrac{\di \xx}{\di x}= e_1 \cos\phi + e_2\sin\phi\qquad
\dfrac{\di \xx}{\di y}=e_1 \cos\phi - e_2 \sin\phi .}
For example, the traveling wave solution
$u=4\arctan\left(\exp(ax+a^{-1}y)\right)$, $a\ne 0$,
gives the classical pseudosphere
(i.e., the surface generated by revolving the tractrix about its asymptote)
when $a^2=1$ and gives Dini's surface when $a^2\ne 1$.

\begin{rem} The relationship between solutions of the \MA/ system for
pseudospherical
surfaces and solutions of sine-Gordon can be described in terms of a double
fibration that is similar to the geometric definition for \BB/ transformations
given in \S\ref{intro}.  For, let $\scM=J^1(\R^2,\R)$ carry the sine-Gordon
system given in \S\ref{SGsec}, and let $\scrU$ be the 5-dimensional quotient
manifold on which $\I$ lives.  Then there is a double fibration
\setlength{\unitlength}{2pt}
\begin{center}
\begin{picture}(40,30)(0,0)
\put(5,5){\makebox(0,0){$\scM$}}
\put(35,5){\makebox(0,0){$\scrU$}}
\put(20,25){\makebox(0,0){$\scM \times \F$}}
\put(16,21){\vector(-3,-4){9}}
\put(24,21){\vector(3,-4){9}}
\put(7,18){\makebox(0,0){$\scriptstyle{\pi_1}$}}
\put(33,18){\makebox(0,0){$\scriptstyle{\pi_2}$}}
\end{picture}
\end{center}
with the obvious projections.  Moreover, \eqref{Aeq} and \eqref{xeq}
are equivalent
to a Pfaffian system $\calS$ of rank six on $\scM\times \F$, with the
property that
if $\Sigma \subset \scM$ is any integral surface of the sine-Gordon system,
then $\calS$ restricts to be Frobenius on $\scN_\Sigma=\pi_1^{-1}(\Sigma)$.
Thus, $\scN_\Sigma$ is
foliated by integral surfaces of $\calS$, each of which projects via $\pi_2$ to
give an integral surface of $\I$.   For this reason, we say that the
pseudospherical
system is an {\em integrable extension} \cite{BG2} of the sine-Gordon system.
However, the double fibration described here does not give a \BT/,
in the sense of \S\ref{intro}, between the two \MA/ systems.
\end{rem}

\subsection*{\BB/ transformations}
The \BB/ transformation for pseudospherical surfaces arises naturally in
a classical context, the study of {\em line congruences} in Euclidean space.
A line congruence is a two-parameter family of lines; associated to the
congruence are two {\em focal surfaces} to which each line in the family
is tangent.  Leaving degeneracies aside, locally the lines give a
1-to-1 correspondence between points on the focal surfaces.

\begin{thm*}[\BB/]
If the distance $\lambda$ between corresponding points on the focal surfaces
and the angle $\psi$ between the surface normal at corresponding
points are both
constant, then the two surfaces have the same constant negative Gauss
curvature,
and we say they are related by a \BT/.
\end{thm*}
(For a proof using moving frames, see \cite{chernterng}.)

For the rest of this section, we normalize the Gauss curvature to be $-1$.
Then $\lambda=\sin\psi$, and we may
regard $\lambda$ or $\psi$ as the parameter in the \BT/.

Let $\F$ and $\Fbar$ be two copies of the frame bundle.  The starting point
of the proof of \BB/'s Theorem is adapting frames $(\xx,e_1,e_2,e_3)$ and
$(\xxbar,\ebar_1,\ebar_2,\ebar_3)$ along the two surfaces so that
$e_1=\ebar_1$ is tangent to the line connecting corresponding points
$\xx$ and $\xxbar$.
Then the graph of the \BT/ is a 6-dimensional submanifold of
$\F\times \Fbar$ on which
\begin{equation}\label{pseudorel}
\begin{aligned}
\bar{\xx} &=\xx+\lambda e_1\\ \ebar_1 &= e_1 \\ \ebar_2 &=
e_2\cos\psi + e_3\sin\psi
\\ \ebar_3 &=e_3\cos\psi  - e_2\sin\psi
\end{aligned}
\end{equation}
for some constant $\lambda = \sin \psi$.

We will regard \eqref{pseudorel} as defining a
7-dimensional submanifold $\scP\subset\F\times \Fbar \times \R$.
Differentiating
the above equations reveals the following relationships between the
canonical forms:
\begin{align*}
\etabar^1 &= \eta^1 + d\lambda, \qquad \etabar^2_3 =\eta^2_3 -d\psi, \\
\begin{pmatrix}\etabar^2 \\ \etabar^3 \end{pmatrix}
        &= \begin{pmatrix}\cos\psi & \sin\psi  \\ -\sin\psi &
\cos\psi\end{pmatrix}
           \begin{pmatrix} \eta^2-\lambda\eta^1_2 \\ \eta^3-\lambda
\eta^1_3 \end{pmatrix} \\
   \begin{pmatrix}\etabar^1_2 \\ \etabar^1_3 \end{pmatrix}
        &= \begin{pmatrix}\cos\psi & \sin\psi  \\ -\sin\psi &
\cos\psi\end{pmatrix}
           \begin{pmatrix} \eta^1_2 \\ \eta^1_3\end{pmatrix}.
\end{align*}

\subsection*{$G$-Structure Invariants}
On $\scP$ we will choose $\theta_1,\theta_2$ to be multiples of
$\eta^3,\etabar^3$ respectively,
and let $\calK=\{\eta^3,\etabar^3,d\psi\}$.
That \eqref{pseudorel} really does define a \BT/ for pseudospherical
surfaces follows from
$$d\etabar^3 \equiv \eta^3_1 \& \eta^3_2 + \eta^1 \& \eta^2 \mood \calK.$$
Moreover, this shows us how we may choose a section of the
$G$-structure $\G$ on $\scP$
associated to the \BT/.  The coframe
$$\begin{aligned}\zeta &= d\psi \\ \theta_1 &= 2\eta^3 \\ \theta_2 &=
2\etabar^3\end{aligned}
\qquad \begin{aligned}\w^1 &=\eta^1-\eta^2_3\\ \w^2 &=\eta^2+\eta^1_3 \\
\w^3 &=\eta^1+\eta^2_3\\ \w^4 &=\eta^2-\eta^1_3\end{aligned}$$
satisfies \eqref{dtzero} with $A_1=-1$ and $A_2=1$.  If we modify
this coframe by setting
\begin{align*}
\w^2 &= \eta^2+\eta^1_3-\dfrac{1+\cos\psi}{\sin\psi}\eta^3 \\
\w^3 &= \etabar^1+\etabar^2_3\\
\w^4 &= \eta^2-\eta^1_3+\dfrac{1-\cos\psi}{\sin\psi}\eta^3
=\etabar^2-\etabar^1_3-\dfrac{1-\cos\psi}{\sin\psi}\eta^3,
\end{align*}
then we get a section of $\G$.  Using this modified coframe, we compute
\begin{alignat*}{3}
\begin{bmatrix} B_1 & C_1 & E_1\\
   B_2 & C_2 & E_2 \end{bmatrix} &=
   \dfrac{1}{4(1-\cos\psi)} \begin{bmatrix} \csc\psi & 0 & 0\\
   0 & -2\sin\psi & 2 \end{bmatrix} \qquad &
\begin{bmatrix} D_1 \\ D_2 \end{bmatrix} &=
\begin{bmatrix} 0 \\ 1+\cos\psi \end{bmatrix} \\[0.1in]
\begin{bmatrix} B_3 & C_3 & E_3 \\
B_4 & C_4 & E_4 \end{bmatrix} &=
\dfrac{-1}{4(1+\cos\psi)}\begin{bmatrix} \csc\psi & 0 & 0\\
   0 & -2\sin\psi & 2\end{bmatrix}\qquad &
\begin{bmatrix} D_3 \\ D_4 \end{bmatrix} &=
\begin{bmatrix} 0 \\ 1-\cos\psi \end{bmatrix}.
\end{alignat*}
Again, $F_1$ and $F_2$ are nonzero, and their ratio is non-constant.
We observe that
\begin{enumerate}
\item{The vector pairs $\{[B_1\ B_2],\ [C_1 \ C_2]\}$ and $\{[B_3 \
B_4],\ [C_3 \ C_4]\}$ are both linearly independent, by contrast with
the sine-Gordon example.}
\item{As in the sine-Gordon example,
the vectors $[D_1\ D_2]$ and $[D_3\ D_4]$ are perpendicular to
their respective $B$-vectors, and the $E$-vectors are dependent
on their respective $C$-vectors.}
\end{enumerate}

\subsection*{Symmetries}
For the pseudospherical surface system on $\F$, the Cauchy
characteristic vector fields correspond to infinitesimal rotation
of the frame within the $e_1e_2$ plane.  Other symmetries generate the
six-dimensional group of Euclidean motions of $\R^3$, and these
vector fields are precisely
the {\em right-invariant} vector fields on $\F$ as a Lie group.

\begin{prop}\label{pseudosymm} These are the only symmetries of the
pseudospherical
system $\I$.
\end{prop}
\begin{proof}Choose a finite-dimensional complement to $\fc_\I \subset \fs_\I$
by requiring that the symmetry vector field satisfy $\vv \intprod \eta^1_2=0$.
We will show that the space of such symmetries is 6-dimensional, and hence
consists only of the projections of the right-invariant vector fields into
the complementary subspace.

Suppose that
$$\vv\intprod \eta^3 = x,\qquad \vv\intprod \eta^i = a_i,\qquad
\vv\intprod \eta^3_i=b_i$$
for $i=1,2$.  Then the requirement that $\Lie_\vv \eta^3 \in \I$ implies
that
\bel{pssymx}{dx \equiv a_i \eta^i_3+b_i \eta^i\mood \eta^3.}
Therefore, let $c$ be the $\eta^3$-component of $dx$.  Differentiation gives
$$0=d(dx)=\alpha_i \& \eta^i_3 + \beta_i \& \eta^i+\gamma \& \eta^3$$
where $\alpha_i = da_i - a_j \eta^j_i$, $\beta_i = db_i - b_j
\eta^j_i-c\eta^3_i$,
$\gamma=dc-b_i \eta^i_3$, and all indices range between 1 and 2.  By the Cartan
Lemma \cite{EDS}, there must be a symmetric matrix $S$ of functions such that
\begin{equation}\label{pssymp}
\transpose{(\alpha_1, \alpha_2, \beta_1, \beta_2, \gamma)}
= S\,
\transpose{(\eta^1_3, \eta^2_3, \eta^1, \eta^2, \eta^3)}
\end{equation}
Let the entries of $S$ be $S_{\mu\nu}=S_{\nu\mu}$ with indices
between 1 and 5.  The
additional symmetry condition is $\Lie_\vv \Psi \equiv 0\mood \I$.
One computes that
$$\Lie_\vv\Psi \equiv 2(S_{13} +S_{24} -c)\eta^1 \&
\eta^2+2(S_{12}+S_{34})\eta^1_3 \& \eta^1
+(S_{11}-S_{44}+x) \eta^1_3 \& \eta^2 -(S_{22}-S_{33}+x)\eta^2_3 \& \eta^1$$
modulo $\eta^3, d\eta^3, \Psi$.
Therefore, we must have $S_{24}=c-S_{13}$,
$S_{34}=-S_{12}$, $S_{44}=S_{11}+x$ and $S_{33}=S_{22}+x$.
We may take \eqref{pssymx},\eqref{pssymp} as defining a Pfaffian
system of rank 6.  Upon calculating the prolongation of this system,
we find that the
only integral manifolds satisfying the independence condition
$\eta^1_3 \& \eta^2_3 \& \eta^1 \& \eta^2 \& \eta^3 \ne 0$ have $c=0$,
and satisfy the rank 6 Frobenius system defined by \eqref{pssymx} and
\begin{align*}
da_1 &= -x\eta^1_3+S_{14}\eta^2-b_1\eta^3-a_2\eta^1_2\\
da_2 &= -x\eta^2_3-S_{14}\eta^1-b_2\eta^3+a_1\eta^1_2\\
db_1 &= -S_{14}\eta^2_3-b_2\eta^1_2\\
db_2 &= S_{14}\eta^1_3+b_1\eta^1_2\\
dS_{14} &=b_1 \eta^2_3 - b_2 \eta^1_3.
\end{align*}
Since any Frobenius system of rank $k$ is locally equivalent to
a system of $k$ first-order ODE, it follows that the space of
such symmetries is parametrized by six constants, and hence is spanned
by the projections of the Euclidean symmetries into the complement of $\fc_\I$.
\end{proof}

\subsection*{Conservation laws}
First we compute the space of conservation laws for the ideal $\calI_1$.
(This will, of course, be isomorphic to the space of conservation
laws for $\calI_2$.)
Suppose that
\[ \Phi = Q\,[\eta^3_1 \& \eta^3_2 + \eta^1 \& \eta^2] + \eta^3 \& \gamma \]
is a closed form in $\calI_1$.  Setting
\[ dQ = Q_1\, \eta^1 + Q_2\, \eta^2 + Q_3\, \eta^3 + Q_4\, \eta^3_1 +
Q_5\, \eta^3_2 \]
and computing
$d\Phi \equiv 0 \mod{\eta^3}$
shows that
\[ \gamma = -Q_5\, \eta^1 + Q_4\, \eta^2 + Q_2\, \eta^3_1 - Q_1\, \eta^3_2. \]
Then the condition $d\Phi = 0$ gives the first-order PDE
\[ Q_3 = 0 \]
and four additional second-order PDE's for $Q$.  The compatibility
conditions for this system lead to five additional equations; the
resulting system is not involutive, so it must be prolonged.
The resulting system is a Frobenius Pfaffian system of rank six;
thus the space of solutions (and hence the space of conservation laws
for each of $\calI_1$ and $\calI_2$) is 6-dimensional.
This is the expected result: the pseudospherical system is variational,
so by N\"oether's theorem there is a one-to-one correspondence between
conservation laws and symmetries.

Next we compute the space of conservation laws for $\calJ$,
using the coframing described above, i.e.
\[ \begin{aligned}
\theta_1 & = 2\eta^3, \\
\w^1 & = \eta^1-\eta^2_3, \\
\w^2 & = \eta^2+\eta^1_3-\dfrac{1+\cos\psi}{\sin\psi}\eta^3,
\end{aligned}\qquad
\begin{aligned}
\theta_2 & = 2\etabar^3, \\
\w^3 & = \etabar^1+\etabar^2_3, \\
\w^4 & = \eta^2-\eta^1_3+\dfrac{1-\cos\psi}{\sin\psi}\eta^3.
\end{aligned} \]
Suppose that
\[ \Phi = \theta_1 \& (P_1\, \w^1 + P_2\, \w^2 + P_3\, \w^3 + P_4\,
\w^4) + \theta_2\& (Q_1\, \w^1 + Q_2\, \w^2 + Q_3\, \w^3 + Q_4\,
\w^4) + R\, \theta_1 \& \theta_2 \]
is a closed form in $\calJ$.  Computing
$d\Phi \equiv 0 \mod{\theta_1, \theta_2}$
shows that
$Q_1 = -P_1$, $Q_2 = -P_2$, $Q_3 = P_3$, and $Q_4 = P_4$.
Next, computing
$d\Phi \equiv 0 \mod{\w^3, \w^4, \theta_1 - \theta_2}$
shows that
\[ R = \frac{(1 + \cos \psi1)\, P_2 + (1 - \cos \psi)\, P_4}{2\sin \psi}. \]
Then the condition $d\Phi=0$ gives a system of 14 first-order PDE's
for the four functions $P_1, P_2, P_3, P_4$.  The compatibility
conditions for this system lead to seven additional equations; the
resulting system is not involutive, so it must be prolonged.
The resulting system is a Frobenius Pfaffian system of rank seven;
thus the space of solutions (and hence the space of conservation laws
for $\calI_1$) is 7-dimensional.

Now we ask the same question as in the sine-Gordon example:  given a
conservation
law for $\calI_1$, is its pullback to $\calJ$ equivalent to the
pullback of a conservation law for $\calI_2$?  By contrast with that case,
the answer here is yes.  We can explain this intuitively using the
duality between conservation laws and symmetries for variational
\MA/ systems.  The six independent
symmetries of $\calI_1$ are all symmetries of $\scB$ as well, unlike in
the sine-Gordon case where the Lie transformation $L_{\mu}$ is not a
  symmetry of $\scB$.  Thus the pullbacks of the spaces of conservation laws
  for $\calI_1$ and $\calI_2$ span the same 6-dimensional subspace of the
  7-dimensional space of conservation laws for $\calJ$.  The
significance of the ``extra"
  conservation law for $\calJ$ is not clear; we hope to explore this
issue in the future.

\section{Timelike CMC surfaces in Lorentzian quadrics}\label{CMCsec}
In this section we study parametric \BB/ transformations for timelike
surfaces of constant mean curvature
in the standard negatively curved Lorentzian space form.  This
example will show
several similarities to the pseudospherical example, but exhibits
values for the
\BB/ invariants which are different from any of the previous examples.
Before computing these invariants, we will briefly review surface
theory in this setting
and construct the \BT/.

Let $\LL^{2,2}$ denote $\R^4$ equipped with an inner product with
$++--$ signature, and
let $\mathbf{H}^{2,1}$ be the quadric hypersurface in $\LL^{2,2}$
defined by $\langle \xx,\xx \rangle =-1$.  Then $\mathbf{H}^{2,1}$
inherits a Lorentzian metric
from $\LL^{2,2}$ and is a homogeneous space form under the action of $O(2,2)$.
Let $\F$ be the bundle of orthogonal
frames $(\xx,e_1,e_2,e_3)$ such that $\xx\in \mathbf{H}^{2,1}$ and
$e_i\in T_\xx \mathbf{H}^{2,1}$, with $\langle e_1,e_1 \rangle =
\langle e_3, e_3 \rangle = 1$ and $\langle e_2,e_2\rangle =-1$.
On $\F$ we define canonical 1-forms $\eta^i, \eta^i_j$ satisfying
$$d\xx = e_i \eta^i, \qquad de_i = \xx \eta^0_j+e_j \eta^j_i,$$
where $\eta^0_1=\eta^1$, $\eta^0_2=-\eta^2$, $\eta^0_3=\eta^3$,
$\eta^1_2=\eta^2_1$,
$\eta^1_3=-\eta^3_1$, $\eta^2_3=\eta^3_2$, and the diagonal
$\eta^i_j$'s are zero.
The structure equations
$$d\eta^i=-\eta^i_j\& \eta^j,\qquad d\eta^i_j=-\eta^i_k \& \eta^k_j
-\eta^i \& \eta^0_j$$
imply that all planes in $\mathbf{H}^{2,1}$ have sectional curvature $-1$
(see \cite{adultoneill} for curvature conventions).

A {\em timelike} surface $S\subset \mathbf{H}^{2,1}$ is one to which
the metric restricts to have $+-$ signature.
A {\em first-order adapted framing} on $S$ is a (local) lift into
$\F$ such that $e_3$ is
normal to the tangent plane of $S$.  Such lifts may be modified by a Lorentzian
rotation of $e_1$ and $e_2$.  Along such a lift $\Sigma \subset \F$
we have $\eta^3=0$,
and hence \eqref{heq} holds
for a symmetric matrix of functions $h_{ij}$.
Computing the change in $h_{ij}$ under rotations shows that
$\det(h_{ij})$ and $h_{11}-h_{22}$
are invariants.  In fact, the Gauss and mean curvature of $S$ are defined by
\begin{align}
d\eta^1_2 = -K\eta^1 \& \eta^2  &=-\eta^1_3 \& \eta^3_2 + \eta^1\&
\eta^2 \label{geqn}
\\ \notag
2H\eta^1 \& \eta^2 &= \eta^3_1 \& \eta^2+\eta^3_2 \& \eta^1.
\end{align}

\subsubsection*{The CMC System}
Suppose $S$ has constant mean curvature $H$.  Then $\Sigma$ is an
integral surface
of the \MA/ system generated by $\eta^3$, $d\eta^3$ and
$$\Psi=\eta^1_3 \& \eta^2 - \eta^2_3\& \eta^1 +2H\eta^1 \& \eta^2.$$
Integral surfaces of this system correspond to solutions of a
certain second-order hyperbolic PDE in local coordinates on $S$ obtained
by choosing first-order adapted frames which diagonalize $h_{ij}$.
For, suppose that $h_{12}=0$,
   $h_{11}=e^{2u}+H$ and $h_{22}=e^{2u}-H$ for some function $u$.
Differentiating \eqref{heq} shows that $e^u \eta^1$ and $e^u \eta^2$
are closed,
and hence equal to $dx$ and $dy$, respectively, for some functions
$x$ and $y$.  Determining
$\eta^1_2$ in terms of $u$ and substituting in the Gauss equation \eqref{geqn}
shows that, as a function of the local coordinates $x$ and $y$, $u$
satisfies the equation
$$u_{xx}-u_{yy}=e^{2u}+(1-H^2)e^{-2u}.$$
This is equivalent to the sinh-Gordon equation when $H^2>1$, the cosh-Gordon
equation when $H^2<1$, and to Liouville's equation when $H^2=1$.
Given a solution $u(x,y)$ of one of these equations defined on an
open set $U\subset \R^2$,
we can construct the corresponding CMC immersion by integrating a
compatible system of
total differential equations like \eqref{Aeq} and \eqref{xeq}.
(In practice, the immersion will only be defined on a sufficiently
small neighbourhood
of a given point in $U$.)

\begin{rem} The correspondence between linear Weingarten surfaces in
Riemannian and
semi-Riemannian space forms and solutions of a handful of `model'
second-order PDE's is set out
in detail in the dissertation of Penn \cite{gabipenn}.
\end{rem}

\subsection*{\BB/ Transformations}
A geometrically-defined \BT/ for timelike CMC surfaces is provided by
the following
analogue of \BB/'s theorem, which is also a generalization of a \BB/
theorem for timelike minimal surfaces in
flat Minkowski space \cite{C01a}.
\begin{thm} Let $U\subset \R^2$ be an open set and
suppose $\xx,\xxbar:U\to \mathbf{H}^{2,1}$ are two timelike immersions with
unit normal vectors $e_3$ and $\ebar_3$ respectively, such that
\begin{list}{(\roman{enumi})}
{\usecounter{enumi}}
\item for each $p\in U$, $\xx(p)$ and $\xxbar(p)$ are the endpoints
of a spacelike geodesic
of fixed length $r>0$;
\item at each end, the tangent vectors to the geodesic orthogonally
project to null vectors in the tangent spaces to the surfaces
$S=\xx(U),\Sbar=\xxbar(U)$;
\item the inner product of the parallel transport of $e_3$
along the geodesic from $\xx$ to $\xxbar$ with the aforementioned
   null vector in $T\Sbar$ is a constant $h \ne 0$.
\end{list}
Then both surfaces have constant mean curvature
$H=\pm\dfrac{2-h}{2\tanh r}.$
\end{thm}
\begin{proof}We may choose first-order adapted frames at
corresponding points on the surface
so that $e_1+e_2+e_3$ and $\ebar_1+\ebar_2+\ebar_3$ are unit tangents to the
geodesic which are related by parallel transport along the geodesic.
It follows that the adapted frames
must be related by equations of the form
\begin{equation}\label{CMCtransform}
\begin{pmatrix}\xxbar \\ \ebar_3 \\ \ebar_1 \\ \ebar_2 \end{pmatrix}
= \begin{pmatrix} c & s & s & s \\ s & (1-h)c & c & (1-h)c \\
                    s & c & (1-k)c & (1-k)c \\
                    -s & (h-1)c &(k-1)c& (h+k-1)c
\end{pmatrix} \begin{pmatrix}\xx \\ e_3 \\ e_1 \\e_2 \end{pmatrix},
\end{equation}
where $c=\cosh r$, $s=\sinh r$, and $c^2 h k =1$.
Computing and comparing the canonical forms shows that
$$0=d\etabar^3=\etabar^1_3 \& \etabar^1 - \etabar^2_3 \& \etabar^2
=\eta^1_3 \& \eta^2 - \eta^2_3\& \eta^1 +\dfrac{2-h}{\tanh r}\eta^1
\& \eta^2.$$
Thus, surface $S$ has the desired mean curvature, and the mean
curvature of $\Sbar$ follows
by symmetry.
(Note, however, that the sign of the mean curvature depends on the orientation
of the normal $e_3$ relative to $e_1$ and $e_2$.  Since the limit of
\eqref{CMCtransform}
as $r\to 0$ interchanges $e_1$ and $e_3$, the sign of the curvature of
the surface $\Sbar$, using normal $\ebar_3$, is {\em opposite} from
that of $S$.)
\end{proof}
Thus, for any fixed value of $H$, there is a one-parameter family of
\BB/ transformations for
timelike CMC surfaces in $\mathbf{H}^{2,1}$.  (Note that we may also
allow $r<0$ as well
as $r>0$ in \eqref{CMCtransform}.)
For the rest of this section, we will restrict to $H=1$, and so we will set
$$h=2(1+\tanh r),\qquad k=\frac12(1-\tanh r).$$

As in \S\ref{pseudosec}, we regard \eqref{CMCtransform} as defining a
7-dimensional
submanifold $\scP\subset\F\times \Fbar \times \R$.  Differentiation shows that
the canonical forms, restricted to $\scP$, satisfy
\begin{align*}
\begin{pmatrix}\etabar^1-\etabar^2 \\ \etabar^2_3-\etabar^1_3 \end{pmatrix}
&= \dfrac1\lambda\begin{pmatrix}1+\lambda & 1-\lambda\\ 1-\lambda &
1+\lambda \end{pmatrix}
\begin{pmatrix}\eta^3-\eta^2 \\ \eta^1_2 +\eta^1_3 \end{pmatrix},
\qquad \lambda=e^{2r}\ne 1, \\
\begin{pmatrix}\etabar^3-\etabar^2 \\ \etabar^1_2+\etabar^1_3 \end{pmatrix}
&=\frac14\begin{pmatrix}1+\lambda & 1-\lambda\\ 1-\lambda & 1+\lambda
\end{pmatrix}
\begin{pmatrix}\eta^1-\eta^2 \\ \eta^2_3-\eta^1_3 \end{pmatrix}, \\
\etabar^1+\etabar^3-\etabar^2 &= \eta^1+\eta^3-\eta^2+dr,\\
\etabar^1_2 + \etabar^1_3 - \etabar^2_3 &=-(\eta^1_2 + \eta^1_3 - \eta^2_3+dr).
\end{align*}

\subsection*{$G$-structure Invariants}
On $\scP$, let $\calK=\{\eta^3,\etabar^3,dr\}$.
The above relations imply that
$$d\etabar^3 \equiv (\eta^1_3 +\eta^1) \& \eta^2 -(\eta^2_3 +\eta^2)
\& \eta^1 \mood \calK$$
while
$$d\eta^3 \equiv (\eta^1_3+\eta^1)\& \eta^1 - (\eta^2_3+\eta^2)\&
\eta^2 \mood \calK.$$
Thus, the coframe
$$\begin{aligned}\zeta&=dr,\\ \theta_1&= 2\eta^3,\\
\theta_2&=2\etabar^3,\end{aligned}
\qquad\begin{aligned} \w^1&=\eta^1-\eta^2,\\ \w^2& =\eta^2_3
+\eta^2+\eta^1_3 + \eta^1, \\
\w^3&=\eta^1+\eta^2, \\ \w^4&=\eta^2_3 +\eta^2-\eta^1_3-\eta^1\end{aligned}$$
satisfies \eqref{dtzero} with $A_1=-1$ and $A_2=1$.  If we modify
this coframe by setting
\begin{align*}
\w^1&=\eta^1-\eta^2+\eta^3\\
\w^3&=\eta^1+\eta^2-\eta^3+dr\\
\w^4&=\eta^2_3 -\eta^1_3-\eta^1+\eta^2+\frac2{\lambda-1}(\eta^3+dr),
\end{align*}
then we get a section of $\G$.  With respect to this coframe, we compute
\begin{alignat*}{3}
\begin{bmatrix} B_1 & C_1 & E_1\\
   B_2 & C_2 & E_2 \end{bmatrix} &=
   \dfrac{1}{4(\lambda-1)} \begin{bmatrix} \lam & 2(\lam-1) & 2\lam \\
   0 & 0 & 0 \end{bmatrix} \qquad &
\begin{bmatrix} D_1 \\ D_2 \end{bmatrix} &=
\begin{bmatrix} 0 \\ 1 \end{bmatrix} \\[0.1in]
\begin{bmatrix} B_3 & C_3 & E_3 \\
B_4 & C_4 & E_4 \end{bmatrix} &=
\dfrac{1}{4(\lam-1)}\begin{bmatrix} \lam & 2(\lam-1) & -2\lam\\[0.1in]
   \dfrac{2\lam}{\lam-1}& -4 & \dfrac{4\lam}{\lam-1}\end{bmatrix}\qquad &
\begin{bmatrix} D_3 \\ D_4 \end{bmatrix} &=
\begin{bmatrix} \dfrac2{\lam-1} \\[0.15in] -1\end{bmatrix}.
\end{alignat*}
As well, $F_1=2/(\lambda-1)$ and $F_2=0$.  In this example, we observe that
\begin{enumerate}
\item{The vectors $\{[B_1\ B_2],\ [C_1 \ C_2]\}$ are linearly dependent,
while $\{[B_3 \ B_4],\ [C_3 \ C_4]\}$ are linearly independent.}
\item{Each $D$-vector is perpendicular to the corresponding $B$-vector.}
\item{The vectors $[E_1\ E_2]$ and $[E_3\ E_4]$ are each linearly dependent on
the corresponding $C$-vectors.}
\end{enumerate}

\subsection*{Symmetries}
A computation similar to that in the proof of Prop. \ref{pseudosymm} yields:
\begin{prop} When $H^2\ne 1$, the only symmetries of the CMC system
for timelike surfaces in $\mathbf{H}^{2,1}$
arise from rigid motions of $\mathbf{H}^{2,1}$ and Lorentzian rotations of the
first-order adapted frames.  When $H^2 = 1$, the space of symmetries depends
(in the sense of Cartan-K\"ahler) on two arbitrary functions of one variable.
\end{prop}

\subsection*{Conservation laws}
First we compute the space of conservation laws for the ideal $\calI_1$.
(This will, of course, be isomorphic to the space of conservation
laws for $\calI_2$.)
Suppose that
\[ \Phi = Q\,[\eta^1_2 \& \eta^2 - \eta^2_3 \& \eta^1 + 2\eta^1 \&
\eta^2] + \eta^3 \& \gamma \]
is a closed form in $\calI_1$.  Setting
\[ dQ = Q_1\, \eta^1 + Q_2\, \eta^2 + Q_3\, \eta^3 + Q_4\, \eta^3_1 +
Q_5\, \eta^3_2 \]
and computing
$d\Phi \equiv 0 \mod{\eta^3}$
shows that
\[ \gamma = (Q_2 - 2Q_5)\, \eta^1 + (Q_1 + 2Q_4)\, \eta^2 + Q_5\,
\eta^3_1 + Q_4\, \eta^3_2. \]
Then the condition $d\Phi = 0$ gives the first-order PDE
\[ Q_3 = 0 \]
and four additional second-order PDE's for $Q$.  The compatibility
conditions for this system lead to three additional equations; the
resulting system is not involutive, so it must be prolonged.
The prolonged system is involutive with last nonvanishing Cartan
character $s_1=2$,
so the space of solutions (and hence the space of conservation laws for ecah of
$\calI_1$ and $\calI_2$) depends on two arbitrary functions of one variable.
This is the expected result: the CMC system is variational, so by
N\"oether's theorem
there is a one-to-one correspondence between conservation laws and symmetries.

Next we compute the space of conservation laws for $\calJ$, using
the coframing described above, i.e.
\[ \begin{aligned}
\theta_1 & = 2\eta^3, \\
\w^1 & = \eta^1-\eta^2+\eta^3, \\
\w^2 & = \eta^2_3 +\eta^2+\eta^1_3 + \eta^1,
\end{aligned}\qquad
\begin{aligned}
\theta_2 & = 2\etabar^3, \\
\w^3 & = \eta^1+\eta^2-\eta^3,  \\
\w^4 & = \eta^2_3 -\eta^1_3-\eta^1+\eta^2+\frac2{\lambda-1}\eta^3.
\end{aligned} \]
Suppose that
\[ \Phi = \theta_1 \& (P_1\, \w^1 + P_2\, \w^2 + P_3\, \w^3 + P_4\,
\w^4) + \theta_2\& (Q_1\, \w^1 + Q_2\, \w^2 + Q_3\, \w^3 + Q_4\,
\w^4) + R\, \theta_1 \& \theta_2 \]
is a closed form in $\calJ$.  Computing
$d\Phi \equiv 0 \mod{\theta_1, \theta_2}$
shows that
$Q_1 = -P_1$, $Q_2 = -P_2$, $Q_3 = P_3$, and $Q_4 = P_4$.
Next, computing
$d\Phi \equiv 0 \mod{\w^3, \w^4, \theta_1 - \theta_2}$
shows that
\[ R = \tfrac{1}{2}(P_1 + P_3) + \frac{1}{\lambda - 1}\, P_4. \]
Then the condition $d\Phi=0$ gives a system of 14 first-order PDE's
for the four functions $P_1, P_2, P_3, P_4$.  The compatibility
conditions for this system lead to four additional equations; the
resulting system is not involutive, so it must be prolonged.
The prolonged system is involutive with last nonvanishing Cartan
character $s_2=1$,
so the space of solutions (and hence the space of conservation laws
for $\calJ$)
depends on one arbitrary function of two variables.

Because the Cartan-K\"ahler analysis does not give explicit expressions
for the conservation laws for $\calI_1$ and $\calI_2$,
it is difficult to determine whether the pullbacks of these spaces to
$\scB$ coincide.
In any case, since the space of conservation laws
for $\calJ$ is strictly larger than those for $\calI_1$ and $\calI_2$,
the pullbacks cannot possibly span the entire space of conservation
laws for $\calJ$.
Again, the significance of the ``extra" conservation laws for $\calJ$
is unclear,
and we hope to explore this issue in the future.

\section{Further Results}\label{concon}
\subsubsection*{Characteristic Systems and their Derived Flags}
We now give some results relating the structure of the
characteristics of the \MA/
systems linked by a parametric \BT/ to the $G$-structure invariants
derived in \S\ref{Gstrsec}.

Recall that $\calC_{11}, \calC_{12}$ denote the characteristic
systems of $\I_1$ and
$\calC_{21}, \calC_{22}$ denote those of $\I_2$, and we have numbered
them so that
$\calC_{11} \equiv \calC_{21}$ and $\calC_{12} \equiv \calC_{22}$
modulo $\calK$.
\begin{thm}\label{intthm} If the vectors $[B_1\ B_2]$, $[E_1 \ E_2]$
are both linearly
dependent on $[C_1 \ C_2]$ at each point of $\G$, but at least two
of these three vectors are nonzero, then each of $\calC_{11}$ and
$\calC_{21}$ contains a
one-dimensional integrable subsystem.  The analogous conditions on
$[B_3\ B_4]$, $[C_3\ C_4]$ and $[E_3 \ E_4]$ imply that $\calC_{12}$
and $\calC_{22}$ each contain
one-dimensional integrable subsystems.
\end{thm}
Before proving this theorem, we will need to explore the implications of
the structure equations \eqref{bigstreq} for the derivatives of the
torsion coefficients
$A_i,\, C_i$, and $D_i$.

Recall that $\calR_1,\calR_2$ denote the Cartan systems of
$\I_1,\I_2$ on $\G$, and
the pullbacks of the characteristic systems to $\G$ include
$\pi_1^*\calC_{11} = \{\theta_1,\w^1,\w^2\}$, $\pi_2^* \calC_{21} =
\{\theta_2,\wt^1,\wt^2\}$.
The fact that the span $\{\theta_2,\wt^1,\wt^2\}$ is the pullback of
a well-defined system
on $\scM_2$ implies that
$$\left.\begin{aligned}d\wt^1 &\equiv 0 \\ d\wt^2 &\equiv 0
\end{aligned}\right\}
\mod \theta_2,\wt^1,\wt^2,\Lambda^2(\calR_2).$$
Substituting in $\wt^1=\w^1+D_2\zeta$ and $\wt^2=\w^2-D_1\zeta$ gives
\begin{equation}\label{Dderivs}
\left.\begin{aligned}
dD_1 &\equiv D_2 \alpha_3 +D_1(\gamma-\alpha_4)
-C_2(D_3\w^3 +D_4\w^4)-E_2\theta_1+D_{10}\theta_2+D_{11}\w^1+D_{12}\w^2\\
dD_2 &\equiv D_1 \alpha_2+D_2(\gamma-\alpha_1)
+C_1(D_3\w^3+D_4\w^4)+E_1\theta_1+D_{20}\theta_2+D_{21}\w^1+D_{22}\w^2
\end{aligned}\right\}\mod \zeta
\end{equation}
for some functions $D_{1j}$ and $D_{2j}$.  (In fact, computing $d(d\theta_2)=0$
and reducing modulo $\theta_2$ and $\w^3$ shows that $D_{12}-D_{21}=F_2$.)
We can derive analogous formulas for $dD_3$ and $dD_4$ by working
with $\calC_{12}$.

Computing $d(d\w^1)=0$ and $d(d\w^2)=0$, and reducing modulo $\w^1,\w^2$ gives
\begin{align}
dC_1 &= C_1(\beta_1+\beta_4-\alpha_1)-C_2\alpha_2+(E_1-C_1F_1)\zeta
+A_2B_1\theta_1 - B_1\theta_2 \notag \\
&  \qquad +C_{11}\w^1+C_{12}\w^2+C_{13}\wt^3 + C_{14}\wt^4 \label{dC1}\\
dC_2 &= C_2(\beta_1+\beta_4-\alpha_4)-C_1\alpha_3+(E_2-C_2F_1)\zeta
+A_2B_2\theta_1 - B_2 \theta_2 \notag \\
&  \qquad +C_{21}\w^2+C_{22}\w^2+C_{23}\wt^3+C_{14}\wt^4 \notag
\end{align}
for some functions $C_{1j}$ and $C_{2j}$.
Then computing $d(d\theta_2) \equiv 0 \mod{\theta_2}$ gives
\begin{equation*}
dA_2 = A_2(\beta_1 + \beta_4 - \alpha_1 - \alpha_4) + A_2 F_2 \zeta
+ C_2 \wt^1 - C_1 \wt^2 + A_{20}\theta_2 + A_{23} \w^3 + A_{24} \w^4
\end{equation*}
for some functions $A_{2j}$.

Recall that the {\em derived system} of a Pfaffian system $\calC$ is
spanned by those
1-forms in $\calC$ whose exterior derivatives are congruent to zero
modulo $\calC$, i.e. are
linear combinations of wedge products, each of which has a form in
$\calC$ as one of its factors.
The derived system, denoted by $\calC^{(1)}$, is equal to $\calC$ if
and only if
$\calC$ is Frobenius.  The derived system of $\calC^{(1)}$ is denoted
by $\calC^{(2)}$, and so on.

It is easy to see from the structure equations \eqref{bigstreq} that
$$\pi_1^*\calC^{(1)}_{11}= \{\w^1-C_1\theta_1,\,\w^2-C_2\theta_2\}.$$
Similarly, using the derivative formulas \eqref{Dderivs}, one computes that
$$\pi_2^*\calC^{(1)}_{21} = \{\wt^1 - (C_1/A_2)\theta_2,\
\wt^2-(C_2/A_2)\theta_2\}.$$
It follows that
$$d(\wt^i-(C_i/A_2)\theta_2)\equiv 0 \mod \wt^1 - (C_1/A_2)\theta_2,\
\wt^2-(C_2/A_2)\theta_2,\ \Lambda^2(\calR_2),$$
and this implies the identities
\bel{Dident}{
\begin{aligned}A_2 D_{10}+C_1 D_{11} + C_2 D_{12} &= D_1 C_{22}-D_2
C_{21}-D_3 C_{24}+D_4 C_{23}+E_2-C_2 F_1\\
A_2 D_{20} + C_1 D_{21} + C_2 D_{22} &=  D_2 C_{11} - D_1 C_{12} +
D_3 C_{14}- D_4 C_{13} - E_1 + C_1 F_1.
\end{aligned}
}
We are now ready to prove Theorem \ref{intthm}.

\begin{proof}  By hypothesis, there exists a sub-bundle $\G'\subset
\G$ on which
$B_1=C_1=E_1 = 0$, and we will restrict all forms and functions to
this sub-bundle.  Then on $\G'$,
$$\pi_1^*\calC^{(1)}_{11}= \{\w^1,\,\w^2-C_2\theta_2\}, \qquad
\pi_2^*\calC^{(1)}_{21} = \{\wt^1,\, \wt^2-(C_2/A_2)\theta_2\}.$$

Computing $d(d\w^1)=0$ on $\G'$ and reducing modulo $\w^1,\w^2$  gives
$$\alpha_2 \& (\theta_1 \& (B_2 \theta_2-E_2\zeta)+C_2\wt^3\&
\wt^4)\equiv 0 \mod \w^1,\w^2.$$
Since $C_2\ne 0$ and at least one of $B_2,E_2$ is nonzero, the 2-form
in parentheses has rank four,
and it follows that $\alpha_2\equiv 0$ modulo $\w^1,\w^2$.  In fact,
substituting into \eqref{dC1}
shows that
\bel{alphaform}{\alpha_2 = \frac{C_{11}\w^1+C_{12}\w^2}{C_2}}
and that $C_{13} = C_{14} = 0$.
It now follows from this and \eqref{bigstreq} that
$$d\w^1 = -\alpha_1 \& \w^1 -\alpha_2 \& \w^2 \equiv 0 \mod \w^1.$$
Thus, $\pi_1^* \calC_{11} = \{\theta_1, \w^1, \w^2\}$ contains the
integrable 1-form $\w^1$.

Furthermore, we compute that
$$d\wt^1\equiv 0 \mod \wt^1,\ \wt^2-(C_2/A_2)\theta_2,$$
which implies that $\wt^1 \in \pi_2^* \calC^{(2)}_{21}$.
(Here, it is necessary to compute keeping the identities
\eqref{Dident} in mind.)
The dimension of this second derived system is upper semi-continuous,
and bounded above by $\dim \calC^{(1)}_{21}=2$.  Hence we may assume that
either $\calC^{(2)}_{21}=\calC^{(1)}_{21}$ on an open set in
$\scM_2$, in which case $\calC^{(1)}_{21}$
is Frobenius and must contain an integrable 1-form, or
$\dim\calC^{(2)}_{21}=1$ on an open set.
In the latter case, then $\pi_2^*\calC^{(2)}_{21} = \{\wt^1\}$, and
so necessarily
$$d\wt^1 \equiv 0 \mod \wt^1, \Lambda^2(\calR_2).$$
However, we compute that
$$d\wt^1 \equiv D_{20}\left(\theta_2 -(A_2/C_2)\wt^2\right) \& \zeta
\mod \wt^1,$$
so it follows in this case that $D_{20}=0$ and $\wt^1$ is integrable.
\end{proof}

\begin{cor} If $[D_1\ D_2]$ is perpendicular to each of
$[B_1\ B_2]$, $[C_1 \ C_2]$ and $[E_1 \ E_2]$, and if at least two of
those three vectors are nonzero,
  then the intersection of $\calC_{11}$ and $\calC_{21}$
contains a one-dimensional integrable subsystem.
Analogous conditions on $[B_3\ B_4]$, $[C_3 \ C_4]$, $[D_3\ D_4]$ and
$[E_3 \ E_4]$
imply that $\calC_{12}$ and $\calC_{22}$ share an integrable 1-form.
\end{cor}
\begin{proof} As in the proof of the previous theorem, we may restrict to
the sub-bundle $\G'$ where $B_1=C_1=E_1$.  Then $D_2=0$, and $\w^1=\wt^1$
is an integrable 1-form common to both systems.
\end{proof}

\subsubsection*{Symmetric \BB/ Transformations}
We will say that a parametric \BT/ $(\scP, \scM_1, \scM_2)$ is {\em
symmetric} if
there is a vector field $\ww$ on $\scP$ such that
\begin{list}{(\alph{enumi})}
{\usecounter{enumi}}
\item $\ww$ is a symmetry of each of the
one-dimensional Pfaffian systems $\{\theta_1\}$, $\{\theta_2\}$ and
$\{\zeta\}$;
\item $\ww$ is transverse to the leaves $\B_\lambda$;
\item $\ww$ has nonzero projections to $\scM_1$ and $\scM_2$.
\end{list}
Although the \BT/s for the sine-Gordon equation and for Goursat's
equation are symmetric,
the \BB/ transformation for pseudospherical surfaces and the
transformation for CMC surfaces with $H^2\ne 1$ are not symmetric,
since the only relevant symmetries of
the underlying systems arise from rigid motions of the ambient space, and
these do not change the \BB/ parameter.  The same argument cannot be applied
to the CMC system with $H=\pm 1$, since in that case the space of  symmetries
(excluding rotations of first-order adapted frames) is not finite-dimensional,
but depends on two functions of one variable.

\begin{prop} The parametric \BT/ for surfaces of constant mean
curvature $H=\pm 1$,
given in \S\ref{CMCsec}, is not symmetric.\end{prop}
\begin{proof}[Sketch] Let $\ww$ be a vector field on $\scP$ and let
$\ww\intprod \w^i=x^i$, $\ww\intprod \w^i_j=y^i_j$ and $\ww\intprod
dr=f$.  Condition (a)
above is equivalent to
\begin{align*}\Lie_\ww du &=g_0du
\\ \Lie_\ww \eta^3 &= g_1 \eta^3
\\ \Lie_\ww \etabar^3 &= g_2 \etabar^3
\end{align*}
for some functions $g_0,g_1,g_2$.  Regard these equations as defining
a rank three Pfaffian
system on $\scP \times \R^{10}$, with $f$ and the $g_i$, $x^i$,
$y^i_j$ as the extra
variables.  The integrability conditions for this Pfaffian system
imply that $f=0$, and
hence there are no vector fields that satisfy both conditions (a) and (b).
\end{proof}

\section*{Concluding Remarks}
There are many issues raised here that merit further exploration.  These include:
\begin{itemize}
\item{In all of the examples computed here, the $E$-vectors and
$C$-vectors are linearly dependent, and the $B$-vectors and $D$-vectors are perpendicular.  Is this true in general?  If
not, what is the significance of these conditions?}
\item{What is the significance of the ``extra" conservation laws for
the system $\calJ$ on $\scB$ appearing in the sine-Gordon and CMC
examples?"}
\item{Is there a condition on the invariants of the $G$-structure
that indicates whether or not a B\"acklund transformation is
symmetric?}
\end{itemize}
We hope to explore these issues in future papers.

\end{document}